\documentclass[a4paper,reqno,11pt]{amsart}
\usepackage[margin=1.2in]{geometry}
\usepackage[onehalfspacing]{setspace}
\usepackage[utf8]{inputenc}
\usepackage[english]{babel}
\usepackage{amsfonts}
\usepackage{amssymb}
\usepackage[reqno]{amsmath}
\usepackage{amsthm}
\usepackage{graphicx}
\usepackage{latexsym}
\usepackage{mathtools}
\usepackage{tikz}
\usepackage{tikz-cd}
\usepackage{setspace}
\usetikzlibrary{cd}
\usepackage[all,cmtip]{xy}
\usepackage[hidelinks]{hyperref}
\usepackage{enumitem}
\usepackage[capitalise,nameinlink,noabbrev]{cleveref}
\usepackage[final]{microtype}
\def\MR#1{} 
\usepackage{todonotes}

\theoremstyle{plain}
\newtheorem{satz}{Satz}[section]
\newtheorem{theorem}[satz]{Theorem}

\newtheorem{lemma}[satz]{Lemma}
\newtheorem{proposition}[satz]{Proposition}
\newtheorem{corollary}[satz]{Corollary}

\theoremstyle{definition}
\newtheorem{definitionx}[satz]{Definition}
\newenvironment{definition}
{\pushQED{\qed}\definitionx}
{\popQED\enddefinitionx}

\newtheorem{examplex}[satz]{Example}
\newenvironment{example}
{\pushQED{\qed}\examplex}
{\popQED\endexamplex}
\newtheorem{remarkx}[satz]{Remark}
\newenvironment{remark}
{\pushQED{\qed}\remarkx}
{\popQED\endremarkx}

\newcommand{\Hom}{\textup{\text{Hom}}}
\newcommand{\End}{\textup{\text{End}}}
\newcommand{\MF}{\textup{\text{MF}}}
\newcommand{\HMF}{\textup{\text{HMF}}}
\newcommand{\Cl}{\textup{\text{Cl}}}

\newcommand{\smod}{\textup{\text{mod}}}
\newcommand{\grmod}{\textup{\text{grmod}}}

\newcommand{\dgCat}{\textup{\text{dgCat}}}
\newcommand{\dgCatD}{\textup{\text{dgCatD}}}
\newcommand{\Sym}{\textup{\text{Sym}}}

\newcommand{\id}{\textup{\text{id}}}

\newcommand{\Ztwo}{\mathbb{Z}\slash 2 \mathbb{Z}}
\newcommand{\Ctwo}{C_2}

\newcommand{\op}{\textup{\text{op}}}
\newcommand{\ev}{\textup{\text{ev}}}
\newcommand{\Aut}{\textup{\text{Aut}}}

\def\pser#1{[\![#1]\!]} 
\newcommand{\ch}{\textup{\text{char}}}
\newcommand{\ptr}{\textup{\text{ptr}}}
\newcommand{\Com}{\textup{\text{Com}}}
\newcommand{\Coh}{\textup{\text{Coh}}}
\newcommand{\Perf}{\textup{\text{Perf}}\,}
\newcommand{\Ind}{\textup{\text{Ind}}\,}

\newcommand{\catC}{\mathcal{C}}
\newcommand{\catD}{\mathcal{D}}
\newcommand{\catE}{\mathcal{E}}

\mathchardef\mhyphen="2D

\title{Matrix factorizations, Reality and Kn\"{o}rrer periodicity}
\author[J.-L. Spellmann]{Jan-Luca Spellmann}
\address{Department of Mathematics and Statistics \\ Utah State University\\
Logan, Utah 84322 \\ USA}
\email{A02370217@usu.edu}

\author[M.\,B. Young]{Matthew B. Young}
\address{Department of Mathematics and Statistics \\ Utah State University\\
Logan, Utah 84322 \\ USA}
\email{matthew.young@usu.edu}

\begin{document}
\date{\today}

\maketitle
	
\begin{abstract}
Motivated by periodicity theorems for Real $K$-theory and Grothendieck--Witt theory and, separately, work of Hori--Walcher on the physics of Landau--Ginzburg orientifolds, we introduce and study categories of Real matrix factorizations. Our main results are generalizations of Kn\"{o}rrer periodicity to categories of Real matrix factorizations. These generalizations are structurally similar to $(1,1)$-periodicity for $KR$-theory and $4$-periodicity for Grothendieck--Witt theory. We use techniques from Real categorical representation theory which allow us to incorporate into our main results equivariance for a finite group and discrete torsion twists.
\end{abstract}
	
\begin{small} \tableofcontents \end{small}
	
\section{Introduction}
\addtocontents{toc}{\protect\setcounter{tocdepth}{1}}
	
The main results of this paper are generalizations of Kn\"{o}rrer periodicity to categories of matrix factorizations which are Real in the sense of Atiyah's Real $K$-theory or symmetric in the sense of its algebraic counterpart, Grothendieck--Witt theory.
	
\subsection{Background and motivation}
Let $R=\mathbb{C}\pser{x_1, \dots, x_n}$ and $w \in R$ a non-zero polynomial without constant term. A matrix factorization of $w$ is a $\Ztwo$-graded finite rank free $R$-module $M$ with an odd $R$-linear endomorphism $d_M$ which satisfies $d_M^2 = w \cdot \id_M$. The $2$-periodic differential graded (dg) category of matrix factorizations $\MF(R,w)$ and its triangulated homotopy category $\HMF(R,w)$ are categorical invariants of the singularity $(R,w)$ introduced by Eisenbud to study the homological algebra of $R \slash (w)$-modules \cite{eisenbud}. Much recent work on matrix factorizations stems from their interpretation as $D$-branes in affine Landau--Ginzburg $B$-models \cite{kapustin2003b}, \cite{walcher2005}, \cite{katzarkov2008}.
    
A fundamental property of matrix factorization categories is Kn\"{o}rrer periodicity.
	
\begin{theorem}[{\cite{knorrer1987}}] \label{thm:knoerreroriginalIntro}
There is a quasi-equivalence of dg categories
\[
\MF(R, w) \xrightarrow[]{\sim} \MF(R\pser{y,z}, w + y^2 + z^2).
\]
\end{theorem}
    
Kn\"{o}rrer periodicity plays an important role in the classification of hypersurface rings of finite maximal Cohen--Macaulay type \cite{knorrer1987}, \cite{leuschke2012} and represents a basic quantum symmetry of Landau--Ginzburg models \cite{herbst2008}. Kn\"{o}rrer periodicity has been generalized in a number of directions. Versions of Kn\"{o}rrer periodicity for global matrix factorizations were proved in \cite{orlov2006}, \cite{shipman2012} and used to study derived categories of projective hypersurfaces \cite{orlov2009}, \cite{ballard2019} and to prove instances of homological projective duality \cite{ballard2018}. Inspired by Homological Mirror Symmetry, an $A$-model counterpart of global Kn\"{o}rrer periodicity was proved for Fukaya--Seidel categories of singular hypersurfaces \cite{jeffs2022}. A different generalization is to matrix factorizations which are equivariant for a group of $\mathbb{C}$-algebra automorphisms of $R$ which preserve $w$. Hirano proved a global version of equivariant Kn\"{o}rrer periodicity and used this to study equivariant derived categories of projective hypersurfaces \cite{hirano2017}. Finally, an $8$-periodic version of Kn\"{o}rrer periodicity for matrix factorizations over $\mathbb{R}$ was proved by Brown used derived Morita theory \cite{brown2016}.

\subsection{Main results}

Recall that the fundamental geometric objects in $KR$-theory are Real vector bundles, that is, complex vector bundles $M$ over a topological space with involution $(X,\sigma)$ with a lift of $\sigma$ to $M$ which is fibrewise $\mathbb{C}$-antilinear \cite{atiyah1966}. In Grothendieck--Witt theory, the algebraic counterpart of $KR$-theory, one instead lifts $\sigma$ to a $\mathbb{C}$-linear isomorphism $M \rightarrow \sigma^* M^{\vee}$ which is symmetric or skew-symmetric \cite{knebusch1977}. According to \cite[\S 7]{dyckerhoff2011}, when $w$ is an isolated hypersurface singularity, $\MF(R,w)$ can be regarded as the dg category of perfect complexes on a hypothetical non-commutative space $\mathcal{X}$ which is dg affine, homologically smooth and proper. Correspondingly, we define two versions of Real matrix factorizations, one antilinear and one contravariant, which can be regarded as the geometric objects of the $KR$-theory and Grothendieck--Witt theory of $\mathcal{X}$. The main results of this paper are Real equivariant generalizations of Kn\"{o}rrer periodicity for both versions of Real matrix factorizations.

To state our results more precisely, let $\Ctwo$ be the multiplicative group $\{1,-1\}$ and $\pi: \hat{G} \rightarrow \Ctwo$ a $\Ctwo$-graded finite group. In the antilinear setting, let $\hat{G}$ act on $R$ by ring automorphisms $\sigma: R \rightarrow R$, $\sigma \in \hat{G}$, which are $\mathbb{C}$-linear if $\pi(\sigma)=1$ and $\mathbb{C}$-antilinear if $\pi(\sigma)=-1$ and such that $w$ is $\hat{G}$-invariant:
\begin{equation} \label{eq:potentialPreserved}
\sigma(w) = w,
\qquad
\sigma \in \hat{G}.
\end{equation}
We stress that while the group $G=\ker \pi$ is a symmetry of $(R,w)$, in the sense of equivariant matrix factorizations discussed above, $\hat{G}$ itself is not because of its antilinear action on $R$. A Real $G$-equivariant structure on a matrix factorization $(M,d_M)$ is the data of coherent isomorphisms $u_{\sigma}: M \rightarrow M$, $\sigma \in \hat{G}$, of graded abelian groups which commute with $d_M$ and satisfy $u_\sigma(rm) = \sigma^{-1}(r)u_\sigma(m)$ for all $r \in R$, $m \in M$. Let $\MF_{\hat{G}}(\mathbb{C}\pser{x_1,\ldots, x_n}, w)$ be the dg category of Real $G$-equivariant matrix factorizations. To state the first form of Real Kn\"{o}rrer periodicity, extend the $\mathbb{C}$-antilinear action of $\hat{G}$ on $R$ to $R\pser{y,z}$ by
\begin{equation}
\label{eq:antilinExt}
\sigma(y)=\pi(\sigma) y, \qquad \sigma(z)=z, 
\qquad \sigma \in \hat{G}.
\end{equation}
Denote by $\Perf(\catC)$ the triangulated hull of a dg category $\catC$.
	
\begin{theorem}[{\cref{thm:RealKnorrer}}] \label{thm:RealKnorrerIntro}
There is a quasi-equivalence of $\mathbb{R}$-linear dg categories
\begin{equation*}
\Perf(\MF_{\hat{G}}(R, w)) \xrightarrow[]{\sim} \Perf(\MF_{\hat{G}}(R\pser{y,z}, w + y^2 + z^2)).
\end{equation*}
\end{theorem}
	
Consider now the contravariant setting, where our approach is motivated by work of Hori--Walcher on the physics of Landau--Ginzburg orientifolds \cite{hori2008}. Let a $\Ctwo$-graded finite group $\hat{G}$ act on $R$ by $\mathbb{C}$-algebra automorphisms such that $w$ is $\pi$-semi-invariant:
\begin{equation} \label{eq:potentialPartiallyPreserved}
\sigma(w) = \pi(\sigma)w,
\qquad
\sigma \in \hat{G}.
\end{equation}
Again, only $G = \ker \pi$ is a symmetry of $(R,w)$. From this data, we construct a duality structure on the dg category $\MF_G(R,w)$ of $G$-equivariant matrix factorizations, that is, a dg functor $\MF_G(R,w)^{\op} \rightarrow \MF_G(R,w)$ with coherence data asserting that it is an involution. Extend the $\mathbb{C}$-linear action of $\hat{G}$ on $R$ to $R\pser{y,z}$ by requiring $G$ to act trivially and
\begin{equation} \label{eq:linExt}
\sigma(y)=-iz, \qquad \sigma(z)=iy, 
\qquad \sigma \in \hat{G} \setminus G.
\end{equation}
The second form of Real Kn\"{o}rrer periodicity is as follows.
		
\begin{theorem}[{\cref{cor:RealKnorrerPeriodicitydgCat,cor:extendedKnorrerPeriodicityCat}}]\label{thm:shiftedKnorrerPeriodicityIntro}
There is a quasi-equivalence of $\mathbb{C}$-linear dg categories with duality
\begin{equation} \label{eq:orthToSympPeriod}
\Perf(\MF_{G}(R,w)) \xrightarrow[]{\sim} \Perf(\MF_G(R\pser{y,z}, w + y^2 + z^2)),
\end{equation}
where the duality structure of the codomain is a shifted and signed version of that of the domain. In particular, there is a quasi-equivalence of dg categories with duality
\begin{equation} \label{eq:fourPeriod}
\Perf(\MF_{G}(R,w)) \xrightarrow[]{\sim} \Perf(\MF_G(R\pser{y_1,z_1,y_2,z_2},w + y_1^2 + z_1^2 + y_2^2 + z_2^2)),
\end{equation}
where both dg categories are given the same duality structure.
\end{theorem}
    
To put \cref{thm:RealKnorrerIntro,thm:shiftedKnorrerPeriodicityIntro} in context, we recall standard periodicities for $KR$-theory and Grothendieck--Witt theory. Atiyah's $(1,1)$-periodicity theorem states that for a compact Hausdorff space $X$ with $\Ctwo$-action, there is an isomorphism
\[
([H]-1) \cdot \; : KR^{\bullet,\bullet}(X) \xrightarrow[]{\sim} KR^{\bullet +1, \bullet +1}(X),
\]
where $[H] -1 \in \widetilde{KR}^0(\mathbb{CP}^1) \simeq KR^{1,1}(\textnormal{pt})$ is the reduced Hopf bundle with $\Ctwo$-action given by complex conjugation \cite{atiyah1966}. As $\Ctwo$-spaces, $\mathbb{CP}^1 \simeq D^{1,1} \slash S^{1,1}$ where $D^{1,1}$ is the closed unit ball in $\mathbb{R} \oplus \mathbb{R}$ with involution $(y,z) \mapsto (-y,z)$ and $S^{1,1}$ is its boundary. That the eigenspaces of the ambient involution with eigenvalue $+1$ and $-1$ each have dimension one is responsible for the $(1,1)$-terminology. Atiyah's $(1,1)$-periodicity unifies many $K$-theoretic periodicities, including Bott's $2$-periodicity for complex $K$-theory and $8$-periodicity for $KO$-theory. In the algebraic setting, associated to a complex affine variety $X$ are its higher Grothendieck--Witt groups $GW^{\bullet}(X)$. These groups are $4$-periodic,
\[
GW^{\bullet}(X) \xrightarrow[]{\sim} GW^{\bullet + 4}(X)
\]
and enjoy natural isomorphisms $GW^0(X) \simeq GW^+(X)$ and $GW^2(X) \simeq GW^-(X)$, where $GW^{\pm}(X)$ is the Grothendieck--Witt group of vector bundles with orthogonal ($+$) or symplectic ($-$) forms.
	
We view \cref{thm:RealKnorrerIntro,thm:shiftedKnorrerPeriodicityIntro} as structural analogues of the $(1,1)$-periodicity and $4$-periodicity theorems, respectively, for matrix factorizations. The $(1,1)$-nature of \cref{thm:RealKnorrerIntro} is apparent from the form of the extension of the $\hat{G}$-action from $R$ to $R\pser{y,z}$; see also the connection with periodicity of categories of Clifford modules below. As for \cref{thm:shiftedKnorrerPeriodicityIntro}, because of the appearance of shifted and signed duality structures, the quasi-equivalence \eqref{eq:orthToSympPeriod} relates bilinear forms of parity $\epsilon \in \Ctwo$ on matrix factorizations of $w$ to bilinear forms of parity $-\epsilon$ on matrix factorizations of $w+y^2+z^2$, and so reflects a shift by two in Grothendieck--Witt theory. The quasi-equivalence \eqref{eq:fourPeriod} then relates bilinear forms of the same parity on matrix factorizations of $w$ and $w + y_1^2 + z_1^2 + y_2^2 + z_2^2$, reflecting the $4$-periodicity of Grothendieck--Witt theory.

\cref{thm:RealKnorrerIntro,thm:shiftedKnorrerPeriodicityIntro} recover a number of known periodicities. Consider first the degenerate case $\hat{G}=\Ctwo$, so that $G$ is trivial. When $\Ctwo$ acts on $R$ by complex conjugation, we recover from \cref{thm:RealKnorrerIntro} the $8$-periodicity of matrix factorizations of non-degenerate quadratic forms over $\mathbb{R}$, as proved by Brown \cite{brown2016}; see \cref{cor:real8Periodicity}. If instead $\Ctwo$ acts on $R$ by a $\mathbb{C}$-algebra involution which negates $w$, the second part of \cref{thm:shiftedKnorrerPeriodicityIntro} is a precise version of the extended Kn\"{o}rrer periodicity discovered by Hori--Walcher \cite{hori2008}. In the degenerate case of trivial $\Ctwo$-grading, so that $\hat{G}=G$, \cref{thm:RealKnorrerIntro,thm:shiftedKnorrerPeriodicityIntro} give an elementary proof of equivariant Kn\"{o}rrer periodicity for affine Landau--Ginzburg models, a special case of Hirano's result \cite{hirano2017}; see \cref{cor:babyHirano}. In Section \ref{sec:quadPot}, motivated by work Buchweitz--Eisenbud--Herzog \cite{buchweitz1987}, we prove in the antilinear setting that the category of Real matrix factorizations of a non-degenerate quadratic form is equivalent to a category of Real graded Clifford modules; see \cref{realbeh}. We use this result in Section \ref{sec:cliffPeriod} to relate \cref{thm:RealKnorrerIntro} to classical periodicities of Clifford modules, including their $(1,1)$-periodicity \cite{atiyah1964}, \cite{lawson1989}, thereby giving a precise sense in which \cref{thm:RealKnorrerIntro} is a $(1,1)$-periodicity theorem.

A second, equally important, context for our results is Landau--Ginzburg orientifolds. Orientifolding is a physical construction which produces an unoriented string theory from an oriented one, in contrast to the standard orbifold construction, which preserves orientability. The mathematics of orientifolds is underdeveloped in comparison to orbifolds. The results of this paper are concrete mathematical consequences of the orientifold construction. In the physics literature, it is well-appreciated that orientifold data in string theory induce on the category of $D$-branes a contravariant involution \cite{diaconescu2007}, \cite{hori2008}. Hori--Walcher identify these involutions for Landau--Ginzburg models and, amongst other things, use them to formulate extended Kn\"{o}rrer periodicity. See \cite{brunner2010} for physical applications of this periodicity. A mathematical approach to some of the ideas of Hori--Walcher was suggested by Bertin--Rosay \cite{bertinrosay} but, as far as the authors are aware, was not pursued. Our approach to \cref{thm:shiftedKnorrerPeriodicityIntro}, and therefore the mathematics of Landau--Ginzburg orientifolds, is different, combining techniques from categorical representation theory, equivariant Grothendieck--Witt theory and the original physical work of Hori--Walcher. While Hori--Walcher and Bertin--Rosay focus on the degenerate case $\hat{G} = \Ctwo$, we work with general $\Ctwo$-graded finite groups $\hat{G}$ and their discrete torsion twists. In this way, we obtain complete results for general affine Landau--Ginzburg orientifolds.
	
Our results suggest two natural $K$-theoretic problems. The first is to relate \cref{thm:RealKnorrerIntro} to the $(1,1)$-periodicity of the $KR$-theory of the Milnor fibre of $w$. For matrix factorizations over $\mathbb{C}$ and $\mathbb{R}$, Brown found such a relation for $K$-theory and $KO$-theory, respectively, under mild assumptions on $w$ \cite{brown2016}. The second is to deduce from \cref{thm:shiftedKnorrerPeriodicityIntro} periodicity for Grothendieck--Witt groups of $G$-equivariant matrix factorizations, as envisioned by Hori--Walcher. At present, the missing ingredient is a $2$-periodic generalization of Schlichting's theorem on the invariance of Grothendieck--Witt groups under dg form equivalence \cite{schlichting2017}.
		
\subsection{Strategy of proof}
We work in the framework of Real $2$-representation theory, as developed by the second author \cite{mbyoung2021b}, \cite{rumyninyoung2021}, in which a $\Ctwo$-graded finite group $\hat{G}$ acts coherently on a $\mathbb{C}$-linear category $\catC$ by functors $\rho(\sigma): \catC \rightarrow \catC$, $\sigma \in \hat{G}$, which are linear and covariant if $\pi(\sigma)=1$ and, depending on the setting, antilinear and covariant or linear and contravariant if $\pi(\sigma)=-1$. In the context of matrix factorizations, a (possibly antilinear) action of $\hat{G}$ on $R$ for which the potential satisfies condition \eqref{eq:potentialPreserved} or \eqref{eq:potentialPartiallyPreserved} defines a Real $2$-representation on $\MF(R,w)$; see \cref{lem:Real2RepMF,lem:contraReal2RepMF}.
    
The quasi-equivalence of \cref{thm:knoerreroriginalIntro} is the dg functor given by tensoring with the matrix factorization
\[
K= \big(
\xymatrixcolsep{2pc}
\xymatrix@C+0.5pc{
\mathbb{C}\pser{y,z} \ar@<+.5ex>[r]^(0.5){y+iz} & \ar@<+.5ex>[l]^(0.5){y-iz} \mathbb{C}\pser{y,z}} 
\big) \in \MF(\mathbb{C}\pser{y,z},y^2+z^2).
\]
Our strategy is to equip $K$, and hence Kn\"{o}rrer's quasi-equivalence by \cref{equivariantfunctor}, with a Real $G$-equivariant structure and then deduce, up taking Real $G$-equivariant objects, \cref{thm:RealKnorrerIntro,thm:shiftedKnorrerPeriodicityIntro}. The details of precisely how this is carried out differ depending on the setting. In both settings, however, a key point is that $K$ admits a Real $G$-equivariant structure only for particular $\hat{G}$-actions on $\mathbb{C}\pser{y,z}$, namely those given by \cref{eq:antilinExt,eq:linExt}; see \cref{prop:RealEquivRankOneMF,prop:RealEquivRankOneMFContra}. This effectively determines the form of our periodicities. There are additional technical difficulties in the contravariant setting. This leads us to prove two general results of independent interest about Real $2$-representations on dg categories (\cref{thm:dualityEquivCat,thm:equivFunctorContra}) which assert that one may take (Real) equivariant objects of Real $2$-representations in stages and that this process is natural.

In general, the triangulated category $\HMF_G(R,w)$ of $G$-equivariant matrix factorizations, and its Real generalizations, is not idempotent complete. For this reason, our techniques lead to statements about $\Perf(\MF_G(R,w))$, whose homotopy category is idempotent complete, instead of $\MF_G(R,w)$ itself. When $w$ is an isolated hypersurface singularity, the case of primary interest for Landau--Ginzburg models, we show in \cref{prop:idemCompHMFG} that $\Perf(\MF_G(R,w))$ is a dg enhancement of the idempotent completion of $\HMF_G(R,w)$. Idempotent completeness is crucial to many applications of matrix factorizations \cite{dyckerhoff2013}, \cite{carqueville2016} and we expect the same to be true equivariantly.

We expect the methods developed in this paper to lead to further mathematical applications of Landau--Ginzburg orientifolds. To mention one, recall that for $w$ an isolated hypersurface singularity, the dg category $\MF(R,w)$ is Calabi--Yau and so defines a $2$-dimensional non-semisimple oriented extended topological field theory \cite{dyckerhoff2012}, \cite{carqueville2016}. We expect that contravariant Real matrix factorizations are the correct framework to lift these theories to unoriented topological field theories. This would confirm the physical predictions of \cite{hori2008}.
        
\subsection*{Acknowledgements}
M.Y. thanks Roland Abuaf, Tobias Dyckerhoff, Eleonore Faber, Gustavo Jasso, Daniel Murfet, Catharina Stroppel and Johannes Walcher for helpful discussions. M.Y. is partially supported by a Simons Foundation Collaboration Grant for Mathematicians (Award ID 853541).
	
\section{Preliminary material}
	
Throughout the paper, $k$ denotes a ground field with $\ch \, k \neq 2$.

\subsection{$2$-periodic dg categories}
    
We recall standard material on differential graded categories, following \cite{keller2006} in the $\mathbb{Z}$-graded setting and \cite[\S 5]{dyckerhoff2011} for $2$-periodic modifications. All categories are assumed to be small by choosing small quasi-equivalent dg categories.

Let $k[u,u^{-1}]$ be the algebra of Laurent polynomials, viewed as a differential $\mathbb{Z}$-graded (dg) algebra with $u$ in degree $2$ and trivial differential. Let $\Com(k)$ be the dg category of complexes of $k$-modules. Finally, let $\Com(k[u,u^{-1}])$ be the symmetric monoidal dg category of dg functors $k[u,u^{-1}] \rightarrow \Com(k)$, where $k[u,u^{-1}]$ is viewed as a dg category with a single object.
    
\begin{definition}
A \emph{$2$-periodic dg category} is a category enriched over $\Com(k[u,u^{-1}])$.
\end{definition}

We often cite results about differential $\mathbb{Z}$-graded categories and apply them in the $2$-periodic setting if their proofs are effectively unchanged. Since we consider only $2$-periodic dg categories in what follows, we abbreviate `$2$-periodic dg' to `dg'.
    
Given a dg category $\catC$, let $Z^0(\catC)$ and $H^0(\catC)$ be the categories with the same objects as $\catC$ and $\Hom_{Z^0(\catC)}(C_1,C_2) = Z^0(\Hom_{\catC}(C_1,C_2))$ and $\Hom_{H^0(\catC)}(C_1,C_2) = H^0(\Hom_{\catC}(C_1,C_2))$.

\begin{definition}
A \emph{dg natural transformation} $\alpha: F \Rightarrow G$ of dg functors $F,G: \catC \rightarrow \catD$ is a family of morphisms $\{\alpha_C \in Z^0(\Hom_{\catD}(F(C),G(C)))\}_{C \in \catC}$ such that $\alpha_{C_2} \circ F(f) = G(f) \circ \alpha_{C_1}$ for each morphism $f: C_1 \rightarrow C_2$. A dg natural transformation whose components are dg isomorphisms, that is, closed degree zero isomorphisms, is called a \emph{dg natural isomorphism}.
\end{definition}
    
We often refer to dg natural transformations simply as natural transformations.
    
Let $\catC$ be a dg category. Its opposite dg category $\catC^{\op}$ has the same objects as $\catC$, morphism complexes $\Hom_{\catC^{\op}}(C_1,C_2) = \Hom_{\catC}(C_2,C_1)$ and composition $f^{\op} \circ g^{\op} = (-1)^{\vert f \vert \vert g \vert} (g \circ f)^{\op}$. Let $\catC^{\op} \mhyphen \smod$ be the dg category of dg functors $\catC^{\op} \rightarrow \Com(k[u,u^{-1}])$. The Yoneda dg functor is $\catC \rightarrow \catC^{\op} \mhyphen \smod$, $C \mapsto \Hom_{\catC}(-,C)$. The pretriangulated hull of $\catC$, denoted $\catC^{\ptr}$, is the smallest dg subcategory of $\catC^{\op} \mhyphen \smod$ which contains the Yoneda image and is closed under shifts and cones of closed morphisms. The triangulated hull of $\catC$, denoted $\Perf(\catC)$, is the full dg subcategory of compact objects of the dg derived category of $\catC^{\op} \mhyphen \smod$.

\begin{definition}
A dg category $\catC$ is called \emph{pretriangulated} (resp. \emph{triangulated}) if the Yoneda dg functor $\catC \rightarrow \catC^{\ptr}$ (resp. $\catC \rightarrow \Perf(\catC)$) is a quasi-equivalence.
\end{definition}
    
Let $F: \catC \rightarrow \catD$ be a dg functor with restriction dg functor $F^* : \catD^{\op} \mhyphen \smod \rightarrow \catC^{\op} \mhyphen \smod$. The derived left adjoint of $F^*$ is the induction dg functor $\Ind F: \Perf(\catC) \rightarrow \Perf(\catD)$. See \cite[\S C]{drinfeld2004} for properties of $\Ind F$.

\subsection{Matrix factorizations} \label{sec:mfCat}
	
We recall background material on matrix factorizations, following \cite{dyckerhoff2011}, \cite{polishchuk2012}.
	
Let $R = k\pser{x_1,\dots, x_n}$ with maximal ideal $\mathfrak{m}$. A non-zero polynomial $w \in \mathfrak{m}$ is called a \emph{potential}. Further restrictions on the potential are common in the literature, depending on the applications in mind. See, for example, \cite[\S 3]{dyckerhoff2011}. One particularly important class of potentials, especially with regards to the physics of affine Landau--Ginzburg models \cite{kapustin2003b}, \cite{walcher2005}, are the isolated hypersurface singularities.

\begin{definition}
\begin{enumerate}[wide,labelwidth=!, labelindent=0pt,label=(\roman*)]
\item The \emph{dg category of matrix factorizations} $\MF(R,w)$ has
\begin{itemize}
\item objects \emph{matrix factorizations} $M=(M, d_M)$, which consists of a finite rank free $\Ztwo$-graded $R$-module $M$ and an odd $R$-linear map $d_M : M \rightarrow M$, the \emph{twisted differential}, which satisfies $d_M^2 = w \cdot \id_M$, and
\item morphism complexes the $\Ztwo$-graded $R$-modules $\Hom_{\MF}(M,N) = \Hom_R(M,N)$ with differential $D$ defined on homogeneous elements by $D(f) = d_N \circ f - (-1)^{\vert f \vert} f \circ d_M$.
\end{itemize}
\item The \emph{homotopy category of matrix factorizations} is $\HMF(R,w) = H^0(\MF(R,w))$.\qedhere
\end{enumerate}
\end{definition}

Let $M \in \MF(R,w)$. After fixing a homogeneous $R$-module decomposition $M=M_0 \oplus M_1$, we can write $d_M = \left( \begin{smallmatrix} 0 & d_M^1 \\ d_M^0 & 0 \end{smallmatrix} \right)$ and depict $M$ as
\[
M= \big(
\xymatrixcolsep{3pc}
\xymatrix@C+0.5pc{
M_0 \ar@<+.5ex>[r]^(0.5){d_M^0} & \ar@<+.5ex>[l]^(0.5){d_M^1} M_1} 
\big).
\]
    
The shift dg functor $\Sigma : \MF(R,w) \rightarrow \MF(R,w)$ is defined on objects and morphisms by $\Sigma(M_0\oplus M_1,d_M)=(M_1 \oplus M_0, -d_M)$ and $\Sigma f = (-1)^{\vert f \vert}f$. Note that $\Sigma \circ \Sigma = \id_{\MF(R,w)}$.
    
Let $w \in R$ and $w^{\prime} \in R^{\prime}$ be potentials. Given $M \in \MF(R,w)$ and $N \in \MF(R',w')$, the \emph{external tensor product} $M \boxtimes N \in \MF(R \otimes_k R', w \otimes 1 + 1 \otimes w')$ is
\[
M \boxtimes N = \big(
\xymatrixcolsep{9pc}
\xymatrix@C+0.5pc{
M_0 \otimes N_0 \oplus M_1 \otimes N_1 \ar@<+.5ex>[r]^(0.5){\left( \begin{smallmatrix} d_M^0 \otimes \id_{N_0} & - \id_{M_1} \otimes d_N^1 \\ \id_{M_0} \otimes d_N^0 & d_M^1 \otimes \id_{N_1} \end{smallmatrix}\right)} & \ar@<+.5ex>[l]^(0.5){\left( \begin{smallmatrix} d_M^1 \otimes \id_{N_0} & \id_{M_0} \otimes d_N^1 \\ -\id_{M_1} \otimes d_N^0 & d_M^0 \otimes \id_{N_1}  \end{smallmatrix}\right)}
M_1 \otimes N_0 \oplus M_0 \otimes N_1 } 
\big).
\]
External tensor product extends to a dg functor 
\begin{equation*}
\boxtimes: \MF(R,w) \otimes_k \MF(R',w') \rightarrow \MF(R \otimes_k R', w \otimes 1 + 1 \otimes w').
\end{equation*}
Under the swap isomorphism $R \otimes_k R^{\prime} \simeq R^{\prime} \otimes_k R$, there is a natural isomorphism
\begin{equation} \label{eq:tensorSym}
M \boxtimes N \simeq N \boxtimes M,
\qquad
m \otimes n \mapsto (-1)^{\vert m\vert \vert n \vert} n \otimes m
\end{equation}
in $\MF(R \otimes_k R^{\prime}, w \otimes 1 +1 \otimes w^{\prime})$. The shift functor and tensor product satisfy
\begin{equation} \label{eq:shiftTensorCompat}
\Sigma M \boxtimes N \simeq M \boxtimes \Sigma N \simeq \Sigma(M \boxtimes N).
\end{equation}
See \cite{yoshino1998} for a detailed discussion of tensor products of matrix factorizations.

The category $\HMF(R,w)$ admits a triangulated structure with shift functor induced by $\Sigma$ \cite[Proposition 3.3]{orlov2004}. If, moreover, $w$ is an isolated hypersurface singularity at the origin, then the dg category $\MF(R,w)$ is triangulated \cite[Lemma 5.6]{dyckerhoff2011}, so that $\HMF(R,w)$ is idempotent complete.
	
A factorization $w=ab$ defines a rank one matrix factorization
\[
\{a,b\} =  \big(
\xymatrixcolsep{3pc}
\xymatrix@C+0.5pc{
R \ar@<+.5ex>[r]^(0.5){a} & \ar@<+.5ex>[l]^(0.5){b} R} 
\big) \in \MF(R,w).
\]
Particularly important to this paper are the cases $w = uv \in k\pser{u,v}$ with associated matrix factorization $\{u,v\}$ and, if $k$ is algebraically closed, $w = y^2 + z^2$ with associated matrix factorization $\{y+iz,y-iz\}$. Knörrer can now be stated as follows.
	
\begin{theorem}[{\cite[Theorem 3.1]{knorrer1987}}]\label{knoerreroriginal}
Let $k$ be algebraically closed with $\ch \, k \neq 2$ and $w \in R= k\pser{x_1,\ldots, x_n}$ a potential. Let $K= \{y+iz, y-iz\}$. The dg functor \begin{equation*}
\mathcal{K} = - \boxtimes K : \MF(R, w) \rightarrow \MF(R\pser{y,z}, w + y^2+z^2)
\end{equation*}
is a quasi-equivalence.
\end{theorem}
	
\subsection{Real $2$-representation theory}\label{sec:Real2Rep}
	
In this section, we establish the representation theoretic framework of the paper, following \cite{mbyoung2021b}.

Denote by $\Ctwo$ the multiplicative group $\{\pm 1\}$. A \emph{$\Ctwo$-graded group} is a group homomorphism $\pi: \hat{G} \rightarrow \Ctwo$. A morphism of $\Ctwo$-graded groups is a group homomorphism which commutes with the $\Ctwo$-gradings. The ungraded group of $\hat{G}$ is $G=\ker\pi$. The terminal $\Ctwo$-graded group is $\id: \Ctwo \rightarrow \Ctwo$, often denoted simply by $C_2$.
		
Let $\hat{G}$ be a $\Ctwo$-graded group with ungraded group $G$. We use two distinct forms of the Real $2$-representation theory of $G$ (with respect to the fixed $\hat{G}$), an antilinear approach in \cref{sec:RealKnorrer,sec:quadPot} and a contravariant approach in \cref{sec:orientifoldMF}. The relevant form will always be clear from the context, so we do not distinguish terminology. When the $\Ctwo$-grading of $\hat{G}$ is trivial, so that $G = \hat{G}$, both forms reduce to earlier notions of categorical representations \cite{deligne1997}, \cite{elagin2014}, \cite{tabauda2018}. We formulate definitions for dg categories, although we sometimes apply them to linear categories, viewed as dg categories concentrated in degree zero.
	
Given a group $G$, denote by $BG$ the locally discrete $2$-category with a single object $\star$ and $1$-morphisms $1\Hom(\star,\star) = G$. A group homomorphism $f: G \rightarrow H$ induces a $2$-functor $B f : B G \rightarrow B H$.
	
\subsubsection{The antilinear approach} \label{sec:antilinApproach}
	
Let $\catC$, $\catD$ be $\mathbb{C}$-linear dg categories. An $\mathbb{R}$-linear dg functor $F:\catC \rightarrow \catD$ is called \emph{antilinear} if its maps on morphism complexes are $\mathbb{C}$-antilinear. Let $\dgCat_{\mathbb{C}/\mathbb{R}}$ be the 2-category whose objects are $\mathbb{C}$-linear dg categories, $1$-morphisms are dg functors which are linear or antilinear and $2$-morphisms are dg natural transformations. There is a $2$-functor $\dgCat_{\mathbb{C}/\mathbb{R}} \rightarrow B\Ctwo$ which records the linearity of the functors.

Fix a $\Ctwo$-graded group $\hat{G}$ with ungraded group $G$.
    
\begin{definition}[{\cite[\S 6.4]{mbyoung2021b}}] \label{def:real2rep}
A \emph{Real $2$-representation} of a $G$ is a $2$-functor $\rho: B\hat{G} \rightarrow \dgCat_{\mathbb{C}/\mathbb{R}}$ over $B\Ctwo$.
\end{definition}
    
Explicitly, a Real $2$-representation $\rho$ is the data of
\begin{enumerate}[label=(\roman*)]
\item a dg category $\rho(\star)=\catC \in \dgCat_{\mathbb{C}/\mathbb{R}}$,
\item $\mathbb{R}$-linear dg functors $\rho(\sigma): \catC \rightarrow \catC$, $\sigma \in \hat{G}$, which are $\mathbb{C}$-linear if $\pi(\sigma)=1$ and $\mathbb{C}$-antilinear if $\pi(\sigma)=-1$,
\item natural isomorphisms $\theta_{\sigma_2,\sigma_1}: \rho(\sigma_2) \circ \rho(\sigma_1) \Rightarrow \rho(\sigma_2 \sigma_1)$, $\sigma_1,\sigma_2 \in \hat{G}$, and
\item \label{cond:unitor} a natural isomorphism $\theta_{e}: \rho(e) \Rightarrow \id_\catC$.
\end{enumerate}
This data is required to satisfy the following coherence conditions:
\begin{enumerate}[label=(\alph*)]
\item \label{cond:coh1} For each $\sigma_1,\sigma_2,\sigma_3 \in \hat{G}$, there is an equality of natural isomorphisms
\begin{equation*}
\theta_{\sigma_3 \sigma_2, \sigma_1} \circ \big(\theta_{\sigma_3, \sigma_2} \circ \id_{\rho(\sigma_1)} \big) = \theta_{\sigma_3, \sigma_2 \sigma_1} \circ \big( \id_{\rho(\sigma_3)} \circ \theta_{\sigma_2, \sigma_1}  \big).
\end{equation*}
\item \label{cond:coh2} For each $\sigma \in \hat{G}$, there are equalities $\theta_{e,\sigma} = \theta_{e} \circ \id_{\rho(\sigma)}$ and $\theta_{\sigma,e} = \id_{\rho(\sigma)} \circ \theta_{e}$ of natural isomorphisms.
\end{enumerate}
		
We often denote a Real $2$-representation by $\rho_{\catC}$ or $(\rho,\theta)$.
	
\begin{definition}
A \emph{Real $G$-equivariant structure} on a dg functor $F: \catC \rightarrow \catD$ between Real $2$-representations of $G$ is a family of natural isomorphisms $\{\eta_{\sigma}: F \circ \rho_{\catC}(\sigma) \Rightarrow \rho_{\catD}(\sigma) \circ F\}_{\sigma \in \hat{G}}$	which makes the diagram
\begin{equation}\label{dia:equivariantfunctor}
\begin{tikzcd}[column sep=5 em,row sep=2.5 em]
F \circ \rho_{\catC}(\sigma_2) \circ \rho_{\catC}(\sigma_1) \arrow[Rightarrow]{r}{\eta_{\sigma_2} \circ \id_{\rho_{\catC}(\sigma_1)}} \arrow[swap,Rightarrow]{d}{F(\theta_{\catC,\sigma_2,\sigma_1})} & \rho_{\catD}(\sigma_2) \circ F \circ \rho_{\catC}(\sigma_1) \arrow[Rightarrow]{r}{\id_{\rho_{\catD}(\sigma_2)} \circ \eta_{\sigma_1}} & \rho_{\catD}(\sigma_2) \circ \rho_{\catD}(\sigma_1) \circ F \arrow[Rightarrow]{d}{\theta_{\catD,\sigma_2,\sigma_1}} \\ 
F \circ \rho_{\catC}(\sigma_2 \sigma_1) \arrow[swap,Rightarrow]{rr}{\eta_{\sigma_2\sigma_1}} & & \rho_{\catD}(\sigma_2\sigma_1) \circ F
\end{tikzcd}
\end{equation}
commute for each $\sigma_1, \sigma_2 \in \hat{G}$.
\end{definition}
	
Denote by $\hat{G} \mhyphen \dgCat_{\mathbb{C} \slash \mathbb{R}}$ the category of Real $2$-representations of $G$ and their Real $G$-equivariant dg functors.

\begin{definition}\label{ex:homotopyFixedCategory}
Let $\catC$ be a Real $2$-representation of $G$. The \emph{homotopy fixed point dg category} $\catC^{\hat{G}}$ has
\begin{itemize}
\item objects \emph{homotopy fixed points} (or \emph{Real $G$-equivariant objects}), which are pairs $(C,u)$ with $C \in \catC$ and $u = \{u_{\sigma}: C \rightarrow \rho(\sigma)(C)\}_{\sigma \in \hat{G}}$ a family of dg isomorphisms such that
\begin{equation}\label{eq:equivariantObject}
u_{\sigma_2 \sigma_1} = \theta_{\sigma_2,\sigma_1} \circ \rho(\sigma_2) (u_{\sigma_1}) \circ u_{\sigma_2}
\end{equation}
for each $\sigma_1,\sigma_2 \in G$, and
	
\item morphisms $f: (C,u) \rightarrow (C',u')$ given by a morphism $f: C \rightarrow C^{\prime}$ in $\catC$ such that
\begin{equation} \label{eq:morphismofhomotopyfixedpoints}
u^{\prime}_{\sigma} \circ f = \rho(\sigma)(f) \circ u_{\sigma}
\end{equation}
for each $\sigma \in \hat{G}$.\qedhere
\end{itemize}
\end{definition}
    
In equation \eqref{eq:equivariantObject} we have written $\theta_{\sigma_2,\sigma_1}$ for what is really its component $\theta_{\sigma_2,\sigma_1,C}$. Similar notational simplifications are made below without comment.
    
Given a Real $G$-equivariant dg functor $(F, \eta): \catC \rightarrow \catD$, let $F^{\hat{G}}: \catC^{\hat{G}} \rightarrow \catD^{\hat{G}}$ be the dg functor defined on objects by $F^{\hat{G}}(X,u)=(F(X), \{\eta_{\sigma} \circ F(u_{\sigma})\}_{\sigma \in \hat{G}})$
and on morphisms by $F^{\hat{G}}(f)=f$. This defines a functor $(-)^{\hat{G}}: \hat{G}\mhyphen \dgCat_{\mathbb{C}\slash \mathbb{R}} \rightarrow \dgCat_{\mathbb{R}}$. See \cite[\S 8]{elagin2014}, \cite[\S 2]{tabauda2018} in the $\mathbb{C}$-linear case, of which the Real case is a direct modification. If the $\Ctwo$-grading of $\hat{G}$ is trivial, then $(-)^{\hat{G}} = (-)^G$ factors through the forgetful functor $\dgCat_{\mathbb{C}} \rightarrow \dgCat_{\mathbb{R}}$.

\begin{proposition} \label{prop:fixedPointFunctorQusiEquiv}
Let $(F,\eta): \catC \rightarrow \catD$ be a Real $G$-equivariant dg functor between Real $2$-representations of a finite group $G$ such that $F$ is a quasi-equivalence.
\begin{enumerate}[wide,labelwidth=!, labelindent=0pt,label=(\roman*)]
\item \label{prop:fixedPointFunctorQusiEquivOne} If $\catC$ and $\catD$ are pretriangulated, then $\Ind F^{\hat{G}}: \Perf(\catC^{\hat{G}}) \rightarrow \Perf(\catD^{\hat{G}})$ is a quasi-equivalence.
\item \label{prop:fixedPointFunctorQusiEquivTwo} If $\catC^{\hat{G}}$ and $\catD^{\hat{G}}$ are triangulated, then $F^{\hat{G}}: \catC^{\hat{G}} \rightarrow \catD^{\hat{G}}$ is a quasi-equivalence.
\end{enumerate}
\end{proposition}

\begin{proof}
\begin{enumerate}[wide,labelwidth=!, labelindent=0pt,label=(\roman*)]
\item Since $\catC$ and $\catD$ are pretriangulated, $H^0(\catC)$ and $H^0(\catD)$ are triangulated and inherit the structure of Real $2$-representations of $G$, with $\hat{G}$ acting by triangle equivalences. It follows from \cite[Corollary 6.10]{elagin2014} that $H^0(\catC)^{\hat{G}}$ and $H^0(\catD)^{\hat{G}}$ are again triangulated. The dg functor $F^{\hat{G}}:\catC^{\hat{G}} \rightarrow \catD^{\hat{G}}$ induces a dg functor $\Ind F^{\hat{G}} :\Perf(\catC^{\hat{G}}) \rightarrow \Perf(\catD^{\hat{G}})$ which fits into the commutative diagram
\begin{equation*}
\begin{tikzcd}[column sep=large,row sep=normal]
H^0(\Perf(\catC^{\hat{G}})) \arrow[r, ""] \arrow[swap, d, "H^0(\Ind F^{\hat{G}})"] & H^0(\catC)^{\hat{G}} \arrow[d, "H^0(F)^{\hat{G}}"]\\
H^0(\Perf(\catD^{\hat{G}})) \arrow[r, ""] & H^0(\catD)^{\hat{G}}.
\end{tikzcd}
\end{equation*}
Since $\catC$ and $\catD$ are pretriangulated, the horizontal arrows are equivalences \cite[Theorem 8.7]{elagin2014}. Since $H^0(F)$ is an equivalence, the same is true of $H^0(F)^{\hat{G}}$ \cite[Proposition 2.20]{sun2019}. Commutativity of the diagram therefore implies that $H^0(\Ind F^{\hat{G}})$ is an equivalence. Since $\Perf(\catC^{\hat{G}})$ and $\Perf(\catD^{\hat{G}})$ are pretriangulated, $\Ind F^{\hat{G}}$ is a quasi-equivalence.

\item If $\catC^{\hat{G}}$ and $\catD^{\hat{G}}$ are triangulated, then the Yoneda arrows in the commutative diagram
\begin{equation*}
\begin{tikzcd}[column sep=large,row sep=normal]
\catC^{\hat{G}} \arrow[r, ""] \arrow[swap, d, "F^{\hat{G}}"] & \Perf(\catC^{\hat{G}}) \arrow[d, "\Ind F^{\hat{G}}"]\\
\catD^{\hat{G}} \arrow[r, ""] & \Perf(\catD^{\hat{G}})
\end{tikzcd}
\end{equation*}
are quasi-equivalences. Part (i) then implies that $F^{\hat{G}}$ is a quasi-equivalence.\qedhere
\end{enumerate}
\end{proof}
    
\begin{remark}\label{rem:goodChar}
\cref{prop:fixedPointFunctorQusiEquiv}, and its proof, apply to $G$-equivariant dg functors between $2$-representations of a finite group $G$ on $k$-linear dg categories, where $k$ is any field with $\ch\, k \nmid \vert G \vert$.
\end{remark}
    
\begin{proposition}\label{equivariantfunctor}
Let $\rho_{\catC}$, $\rho_{\catD}$, $\rho_\catE$ be Real $2$-representations of $G$ and $(F, \eta): \catC \times \catD \rightarrow \catE$ a Real $G$-equivariant dg functor. Each object $(D,u) \in \catD^{\hat{G}}$ induces a Real $G$-equivariant structure on the dg functor $F(-,D)$.
\end{proposition}
\begin{proof}
Define $\eta^D = \{
\eta^D_{\sigma}: F(-,D) \circ \rho_{\catC}(\sigma) \Rightarrow \rho_\catE(\sigma) \circ F(-,D)
\}_{\sigma \in \hat{G}}$ so that the component of $\eta^D_{\sigma}$ at $C \in \catC$ is 
$$\eta^D_{\sigma,C}: F(\rho_\catC(\sigma)(C),D) \xrightarrow{F(\id_C,u_{\sigma})} F(\rho_\catC(\sigma)(C),\rho_\catD(\sigma)(D)) \xrightarrow{\eta_{\sigma}} \rho_\catE(\sigma)(F(C,D)).$$
Naturality of $\eta_{\sigma}$ implies that of $\eta_{\sigma}^D$. That $\eta^D$ satisfies the coherence condition \eqref{dia:equivariantfunctor} follows from the fact that $\eta$ satisfies the same conditions and $u$ satisfies condition \eqref{eq:equivariantObject}.
\end{proof}
    
\subsubsection{The contravariant approach}  \label{sec:contravariantApproach}
	
Let $\dgCat_k$ be the $2$-category of $k$-linear dg categories, dg functors and dg natural transformations. Taking the opposite of a dg category extends to a duality involution of $\dgCat_k$. In particular, if $\theta: F \Rightarrow G$ is a natural transformation, then its opposite is $\theta^{\op} : G^{\op} \Rightarrow F^{\op}$. Given $\epsilon \in \Ctwo$ and $\catC \in \dgCat_k$, write 
\[
{^{\epsilon}}\catC
=
\begin{cases}
\catC & \mbox{if } \epsilon =1, \\
\catC^{\op} & \mbox{if } \epsilon = -1
\end{cases}
\]
with similar notation for dg functors and natural transformations. 

For a direct contravariant analogue of \cref{def:real2rep}, see \cite[\S 4.4]{mbyoung2021b}. We instead use an explicit unpacking of this definition. Fix a $\Ctwo$-graded group $\hat{G}$ with ungraded group $G$.

\begin{definition} \label{def:real2repContra}
A \emph{Real $2$-representation} of $G$ on a dg category $\catC$ is the data of
\begin{enumerate}[label=(\roman*)]
\item dg functors $\rho(\sigma): {^{\pi(\sigma)}}\catC \rightarrow \catC$, $\sigma \in \hat{G}$,
\item natural isomorphisms $\theta_{\sigma_2,\sigma_1}: \rho(\sigma_2) \circ {^{\pi(\sigma_2)}}\rho(\sigma_1) \Rightarrow \rho(\sigma_2 \sigma_1)$, $\sigma_1,\sigma_2 \in \hat{G}$, and
\item a natural isomorphism $\theta_{e}: \rho(e) \Rightarrow \id_\catC$.
\end{enumerate}
This data is required to satisfy the following coherence conditions:
\begin{enumerate}[label=(\alph*)]
\item For each $\sigma_1, \sigma_2, \sigma_3 \in \hat{G}$, there is an equality of natural isomorphisms
\begin{equation} \label{eq:theta2CocycleContra}
\theta_{\sigma_3 \sigma_2,\sigma_1} \circ \big(\theta_{\sigma_3, \sigma_2} \circ \id_{{^{\pi(\sigma_3 \sigma_2)}}\rho(\sigma_1)} \big) = \theta_{\sigma_3,\sigma_2 \sigma_1} \circ \big( \id_{\rho(\sigma_3)} \circ {^{\pi(\sigma_3)}}\theta^{\pi(\sigma_3)}_{\sigma_2,\sigma_1}  \big).
\end{equation}
\item For each $\sigma \in \hat{G}$, there are equalities $\theta_{e,\sigma} = \theta_{e} \circ \id_{\rho(\sigma)}$ and $\theta_{\sigma,1} = \id_{\rho(\sigma)} \circ {^{\pi(\sigma)}}\theta_{e}$ of natural isomorphisms.\qedhere
\end{enumerate}
\end{definition}
    
A \emph{Real $G$-equivariant structure} on a dg functor $F: \catC \rightarrow \catD$ between Real $2$-representations is a family $\{\eta_{\sigma}: F \circ \rho_{\catC}(\sigma) \Rightarrow \rho_{\catD}(\sigma) \circ {^{\pi(\sigma)}}F\}_{\sigma \in \hat{G}}$ of natural isomorphisms such that

\begin{equation} \label{eq:equivariantFunctorContra}
\theta_{\catD,\sigma_2,\sigma_1} \circ \left( \id_{\rho_{\catD}(\sigma_2)} \circ {^{\pi(\sigma_2)}}\eta_{\sigma_1}^{\pi(\sigma_2)} \right) \circ \left(\eta_{\sigma_2} \circ \id_{{^{\pi(\sigma_2)}}\rho_{\catC}(\sigma_1)}\right)
=	
\eta_{\sigma_2 \sigma_1} \circ F(\theta_{\catC,\sigma_2, \sigma_1})
\end{equation}
for each $\sigma_1,\sigma_2 \in \hat{G}$. A \emph{homotopy fixed point} of $\catC$ is a pair $(C,u)$ consisting of an object $C \in \catC$ and dg isomorphisms $u = \{u_{\sigma}: C \rightarrow \rho(\sigma)(C)\}_{\sigma \in \hat{G}}$ such that
\begin{equation*}
u_{\sigma_2 \sigma_1} = \theta_{\sigma_2, \sigma_1} \circ \rho(\sigma_2)(u^{\pi(\sigma_2)}_{\sigma_1}) \circ u_{\sigma_2} \end{equation*}
for each $\sigma_1, \sigma_2 \in \hat{G}$.
    
\begin{proposition}\label{prop:Real2RepOnPerf}
A Real $2$-representation of $G$ on a dg category $\catC$ induces a Real $2$-representation of $G$ on $\Perf(\catC)$ which makes the Yoneda dg functor $\catC \rightarrow \Perf(\catC)$ Real $G$-equivariant. Moreover, a Real $G$-equivariant dg functor $(F,\eta): \catC \rightarrow \catD$ lifts to a Real $G$-equivariant dg functor $(\Ind F, \Ind \eta): \Perf(\catC) \rightarrow \Perf(\catD)$.
\end{proposition}

\begin{proof}
Let $(\rho, \theta)$ be a Real $2$-representation of $G$ on $\catC$. Define a Real $2$-representation $(\rho^{\prime}, \theta^{\prime})$ on $\Perf(\catC)$ as follows. Set $\rho^{\prime}(\sigma) = \Ind \rho(\sigma)$, $\sigma \in \hat{G}$, where, for $\sigma \in \hat{G} \setminus G$ we implicitly use $\Perf(\catC^{\op}) \simeq \Perf(\catC)^{\op}$. The natural isomorphism $\theta_{\sigma_2,\sigma_1}$ induces a natural isomorphism of pullback dg functors
and hence also of their derived left adjoints,
\[
\theta^{\prime}_{\sigma_2,\sigma_1} : \Ind \rho(\sigma_2) \circ {^{\pi(\sigma_2)}}\Ind \rho(\sigma_1) \Rightarrow \Ind \rho(\sigma_2 \sigma_1).
\]
Coherence of $\theta^{\prime}$ follows from that of $\theta$. The identity natural transformations define a Real $G$-equivariant structure on the Yoneda dg functor; see \cite[\S C.10]{drinfeld2004}. The proof of the final statement is similar and so is omitted.
\end{proof}

\subsubsection{Twists by $2$-cocycles} \label{sec:twistCoherenceData}

There is a simple way to construct new Real $2$-representations of $G$ from a given Real $2$-representation. To treat both approaches simultaneously, let $k^{\times}_{\pi}$ be the $\hat{G}$-module
\begin{itemize}
\item Antilinear approach: $\mathbb{C}^{\times}$ with $\hat{G}$ acting through $\pi$ by complex conjugation;
\item Contravariant approach: $k^{\times}$ with $\hat{G}$ acting through $\pi$ by inversion.
\end{itemize}
Let $C^{\bullet}(\hat{G}; k^{\times}_{\pi})$ be the complex of normalized group cochains on $\hat{G}$ with coefficients in $k^{\times}_{\pi}$ and $Z^{\bullet}(\hat{G};k^{\times}_{\pi})$ and $H^{\bullet}(\hat{G}; k^{\times}_{\pi})$ the groups of cocycles and cohomology classes, respectively.

Given $\hat{\mu} \in Z^2(\hat{G}; k^{\times}_{\pi})$ and $(\rho, \theta)$ a Real $2$-representation of $G$ (in either of the two approaches), define a Real $2$-representation $(\rho^{\prime},\theta^{\prime})$ by $\rho^{\prime} = \rho$ and $\theta^{\prime}_{\sigma_2, \sigma_1} = \hat{\mu}([\sigma_2 \vert \sigma_1]) \theta_{\sigma_2,\sigma_1}$, $\sigma_2,\sigma_1 \in \hat{G}$. Up to equivalence, $(\rho^{\prime}, \theta^{\prime})$ depends on $\hat{\mu}$ through $[\hat{\mu}] \in H^2(\hat{G}; k^{\times}_{\pi})$.

The terminal $\Ctwo$-graded group satisfies $H^2(\Ctwo; k^{\times}_{\pi}) \simeq \Ctwo$. A representative of the non-trivial class is $\hat{c}([\sigma_2 \vert \sigma_1])=(-1)^{\frac{\sigma_2-1}{2}\frac{\sigma_1-1}{2}}$. For any $\Ctwo$-graded finite group $\pi: \hat{G} \rightarrow \Ctwo$, pullback along the $\Ctwo$-grading defines $\pi^*: H^2(\Ctwo; k^{\times}_{\pi}) \rightarrow H^2(\hat{G}; k^{\times}_{\pi})$. In particular, associated to any Real $2$-representation $(\rho,\theta)$ of $G$ is a second Real $2$-representation $(\rho, \theta_-)$, where $\theta_-=\pi^*\hat{c} \cdot \theta$.

\section{Kn\"{o}rrer periodicity for Real equivariant matrix factorizations} \label{sec:RealKnorrer}
	
In this section, $k=\mathbb{C}$ and Real $2$-representations are in the antilinear approach of \cref{sec:antilinApproach}. We introduce Real matrix factorizations, in the antilinear setting, and prove the first instance of Real Kn\"{o}rrer periodicity.
	
\subsection{Equivariant module categories} \label{sec:RealEquivModCat}
    
Let $R$ be a $\mathbb{C}$-algebra. The group $\Aut_{\mathbb{C} \slash \mathbb{R}}(R)$ of \emph{generalized algebra automorphisms} of $R$ is the group of ring automorphisms $R \rightarrow R$ which are either $\mathbb{C}$-linear or $\mathbb{C}$-antilinear. The map $\Aut_{\mathbb{C} \slash \mathbb{R}}(R) \rightarrow \Ctwo$ which records linearity makes $\Aut_{\mathbb{C} \slash \mathbb{R}}(R)$ into a $\Ctwo$-graded group with ungraded group $\Aut_{\mathbb{C}}(R)$. Given $\sigma \in \Aut_{\mathbb{C} \slash \mathbb{R}}(R)$, let $(-)^\sigma: R\mhyphen\smod \rightarrow R\mhyphen\smod$ be the functor given on a finitely generated $R$-module by $M \mapsto M^\sigma$, where $M^{\sigma} = M$ as abelian groups with $R$-module structure $r \cdot m = \sigma^{-1}(r) m$, and on an $R$-module homomorphism $f: M \rightarrow N$ by $f \mapsto f^\sigma$, where $f^\sigma(m) = f(m)$, $m \in M^{\sigma}$.

\begin{lemma}\label{2rep}
The functors $\{(-)^\sigma\}_{\sigma \in \Aut_{\mathbb{C} \slash \mathbb{R}}(R)}$ extend to a Real $2$-representation of $\Aut_{\mathbb{C}}(R)$ on $R\mhyphen\smod$. 
\end{lemma}
\begin{proof}
Let $\sigma_1,\sigma_2 \in \Aut_{\mathbb{C} \slash \mathbb{R}}(R)$. The $R$-module structure of $(M^\sigma_1)^{\sigma_2}$ is $r \cdot m = (\sigma_2 \sigma_1)^{-1}(r)m$.
Identity maps of underlying abelian groups therefore define the components of the required $2$-isomorphisms $\{\theta_{\sigma_2,\sigma_1}\}_{\sigma_2,\sigma_1 \in \Aut_{\mathbb{C} \slash \mathbb{R}}(R)}$ and $\theta_e$.
\end{proof}
    
Let $\hat{G}$ be a $\Ctwo$-graded finite group. A $\Ctwo$-graded group homomorphism $f: \hat{G} \rightarrow \Aut_{\mathbb{C} \slash \mathbb{R}}(R)$ is called an \emph{action of $\hat{G}$ on $R$ by generalized algebra automorphisms}. Pullback of the Real $2$-representation of \cref{2rep} along $Bf: B \hat{G} \rightarrow B \Aut_{\mathbb{C} \slash \mathbb{R}}(R)$ defines a Real $2$-representation of $G$ on $R \mhyphen \smod$. The \emph{category of Real $G$-equivariant $R$-modules} is $R\mhyphen\smod^{\hat{G}}$. An object $(M, u) \in R\mhyphen\smod^{\hat{G}}$ is an $R$-module $M$ and a family of $R$-module isomorphisms $\{u_{\sigma}: M \rightarrow M^{\sigma}\}_{\sigma \in \hat{G}}$ which satisfy $u_{\sigma_2 \sigma_1} = (u_{\sigma_1})^{\sigma_2} \circ u_{\sigma_2}$.
More generally, we can twist this Real $2$-representation by $\hat{\mu} \in Z^2(\hat{G}; \mathbb{C}^{\times}_{\pi})$, as in \cref{sec:twistCoherenceData}, in which case an object $(M,u) \in R \mhyphen \smod^{\hat{G}}$ is a $\hat{\mu}$-twisted Real $G$-equivariant $R$-module: the family of $R$-module isomorphisms $\{u_{\sigma}\}_{\sigma \in \hat{G}}$ satisfy $u_{\sigma_2 \sigma_1} = \hat{\mu}([\sigma_2 \vert \sigma_1])(u_{\sigma_1})^{\sigma_2} \circ u_{\sigma_2}$.
	
The case of primary interest for this paper is $R=\mathbb{C}\pser{x_1,\ldots,x_n}$. Using Reynolds operators, we may assume that $\hat{G}$ acts on $R$ by degree preserving generalized algebra automorphisms. In more invariant terms, there exists a finite dimensional Real representation $V$ of $G$, in the sense of \cite[\S 3.1]{mbyoung2021b}, so that $\hat{G}$ acts on $V$ linearly or antilinearly according to the $\Ctwo$-grading, such that $R$ is Real $G$-equivariantly isomorphic to the completed symmetric algebra $\widehat{\Sym\, V^{\vee}}$. We henceforth work in such coordinates.

\subsection{Equivariant matrix factorization categories} \label{sec:RealMFCat}
    
Let $R=\mathbb{C}\pser{x_1, \dots, x_n}$ with potential $w \in \mathfrak{m}$. For each $\sigma \in \Aut_{\mathbb{C} \slash \mathbb{R}}(R)$, let $(-)^\sigma: \MF(R,w) \rightarrow \MF(R,\sigma(w))$ be the dg functor defined on objects and morphisms by $(M,d_M) \mapsto (M^\sigma, d_M^\sigma)$ and  $f \mapsto f^{\sigma}$, respectively.

Let $\hat{G}$ be a $\Ctwo$-graded finite group acting on $R$ by generalized algebra automorphisms which preserve the potential: $\sigma(w) = w$, $\sigma \in \hat{G}$. Each $(-)^{\sigma}$ is then a dg endofunctor of $\MF(R,w)$ and we obtain the following result.
	
\begin{lemma}\label{lem:Real2RepMF}
The dg functors $\{(-)^{\sigma}\}_{\sigma \in \hat{G}}$ extend to a Real $2$-representation of $G$ on $\MF(R,w)$.
\end{lemma}
	
\begin{definition}\label{def:RealMF}
\begin{enumerate}[wide,labelwidth=!, labelindent=0pt,label=(\roman*)]
\item The \emph{category of Real $G$-equivariant matrix factorizations} $\MF_{\hat{G}}(R,w)$ is the homotopy fixed point category $\MF(R,w)^{\hat{G}}$.
\item The \emph{homotopy category of Real $G$-equivariant matrix factorizations} $\HMF_{\hat{G}}(R,w)$ is $H^0(\MF_{\hat{G}}(R,w))$.\qedhere 
\end{enumerate}
\end{definition}

\begin{example}
\label{ex:antilinInvtPot}
\begin{enumerate}[wide,labelwidth=!, labelindent=0pt,label=(\roman*)]
\item As a degenerate case, let a finite group $G$ (with trivial $\Ctwo$-grading) act on $\mathbb{C}\pser{x_1,\dots,x_n}$ by $\mathbb{C}$-algebra automorphisms. Any $G$-invariant potential $w \in \mathbb{C}\pser{x_1,\dots,x_n}$ is then invariant in the above sense. In this example, \cref{def:RealMF} reduces to the $\mathbb{C}$-linear dg category $\MF_G(R,w)$ of $G$-equivariant matrix factorizations, a model for the the category of $D$-branes in the Landau--Ginzburg orbifold associated to the symmetry group $G$ of $(R,w)$ \cite[\S 2.2]{ashok2004}.

\item \label{ex:antilinInvtPotOverR} Let a finite group $G$ act on $\mathbb{R}\pser{x_1,\dots, x_n}$ by $\mathbb{R}$-algebra automorphisms and let $w \in \mathbb{R}\pser{x_1,\dots, x_n}$ be a $G$-invariant potential. Then $w$ is also invariant with respect to the generalized action of $\hat{G} = G \times \Ctwo$ on $\mathbb{R}\pser{x_1,\dots, x_n} \otimes_{\mathbb{R}} \mathbb{C} \simeq \mathbb{C}\pser{x_1,\dots, x_n}$ in which $\Ctwo$ acts on $\mathbb{C}$ by complex conjugation. In particular, when $G$ is the trivial group, taking $\Ctwo$-fixed points defines a dg isomorphism
\[
\MF_{\Ctwo}(\mathbb{C}\pser{x_1, \dots, x_n},w) \xrightarrow[]{\sim} \MF(\mathbb{R}\pser{x_1, \dots, x_n}, w).
\]
In this way, matrix factorizations over $\mathbb{R}$ appear as a special case of Real equivariant matrix factorizations.

\item Let $w = x^m + y^m \in \mathbb{C} \pser{x,y}$, $m \in \mathbb{Z}_{\geq 2}$. Let $\hat{G} = D_{2m}$ be the dihedral group of order $2m$ with $\pi$ the projection with kernel $G$ the cyclic group of order $m$. Define a generalized $\hat{G}$-action on $\mathbb{C} \pser{x,y}$ by letting a generator of $G$ act by $x \mapsto \zeta_m x$ and $y \mapsto \zeta_m y$, with $\zeta_m$ a primitive $m\textsuperscript{th}$ root of unity, and a fixed element of the non-identity coset of $\hat{G}$ act by complex conjugation. Then $w$ is invariant.

\item Let $w \in \mathbb{C}\pser{x_1,\dots, x_n}$ be homogeneous of degree divisible by $4$. Then $w$ is invariant with respect to the generalized action of the terminal $\Ctwo$-graded group determined by $x_j \mapsto i x_j$, $j=1, \dots, n$.\qedhere 
\end{enumerate}
\end{example}
	
\begin{lemma} \label{lem:extTensorEquiv}
\begin{enumerate}[wide,labelwidth=!, labelindent=0pt,label=(\roman*)]
\item \label{lem:extTensorEquivOne} For each $\sigma \in \Aut_{\mathbb{C} \slash \mathbb{R}}(R)$ and $\sigma^{\prime} \in \Aut_{\mathbb{C} \slash \mathbb{R}}(R')$, the diagram
\begin{equation*}
\begin{tikzcd}[column sep=normal,row sep=large]
\MF(R,w) \otimes_{\mathbb{C}} \MF(R',w') \arrow{r}{\boxtimes} \arrow{d}[swap]{(-)^\sigma \otimes (-)^{\sigma^{\prime}}} & \MF(R \otimes_{\mathbb{C}} R', w \otimes 1 + 1 \otimes w') \arrow{d}{(-)^{\sigma \otimes \sigma^{\prime}}} \\
\MF(R,\sigma(w)) \otimes_{\mathbb{C}} \MF(R',\sigma^{\prime}(w')) \arrow{r}[swap]{\boxtimes} & \MF(R \otimes_{\mathbb{C}} R', \sigma(w) \otimes 1 + 1 \otimes \sigma^{\prime}(w'))
\end{tikzcd}
\end{equation*}
commutes up to natural isomorphism.
		
\item With respect to the Real $2$-representations of $G$ of \cref{lem:Real2RepMF}, the dg functor
$$\boxtimes: \MF(R,w) \otimes_{\mathbb{C}} \MF(R',w') \rightarrow \MF(R \otimes_{\mathbb{C}} R', w \otimes 1 + 1 \otimes w')$$
admits a Real $G$-equivariant structure.
\end{enumerate}
\end{lemma}
\begin{proof}
\begin{enumerate}[wide,labelwidth=!, labelindent=0pt,label=(\roman*)]
\item Let $M \in \MF(R,w)$ and $N \in \MF(R',w')$. Direct sums of identity maps of underlying abelian groups give isomorphisms
\begin{equation*}
(M_0 \otimes_\mathbb{C} N_0)^{\sigma \otimes \sigma^{\prime}} \oplus (M_1 \otimes_\mathbb{C} N_1)^{\sigma \otimes \sigma^{\prime}} \rightarrow (M_0^\sigma \otimes_\mathbb{C} N^{\sigma^{\prime}}_0) \oplus (M^\sigma_1 \otimes_\mathbb{C} N^{\sigma^{\prime}}_1)
\end{equation*}
and
\begin{equation*}
(M_1 \otimes_\mathbb{C} N_0)^{\sigma \otimes \sigma^{\prime}} \oplus (M_0 \otimes_\mathbb{C} N_1)^{\sigma \otimes \sigma^{\prime}} \rightarrow (M_1^\sigma \otimes_\mathbb{C} N^{\sigma^{\prime}}_0) \oplus (M^\sigma_0 \otimes_\mathbb{C} N^{\sigma^{\prime}}_1).
\end{equation*}
These are the components of an isomorphism $(M \boxtimes N)^{\sigma \otimes \sigma^{\prime}} \rightarrow M^{\sigma} \boxtimes N^{\sigma^{\prime}}$ that fulfill the coherence condition \eqref{dia:equivariantfunctor}.
\item The identity natural isomorphisms from part \ref{lem:extTensorEquivOne} are the data for the required Real $G$-equivariant structure on $\boxtimes$.\qedhere
\end{enumerate}
\end{proof}
    
\begin{proposition}\label{prop:idemCompHMFG}
Let $w$ be an isolated hypersurface singularity. The canonical functor $\HMF_{\hat{G}}(R,w) \rightarrow \HMF(R,w)^{\hat{G}}$ is an idempotent completion. In particular, $\Perf(\MF_{\hat{G}}(R,w))$ is a dg enhancement of $\HMF(R,w)^{\hat{G}}$.
\end{proposition}

\begin{proof}
The Yoneda dg functor $\MF_{\hat{G}}(R,w) \rightarrow \Perf(\MF_{\hat{G}}(R,w))$ induces on homotopy categories an idempotent completion $\HMF_{\hat{G}}(R,w) \rightarrow H^0(\Perf(\MF_{\hat{G}}(R,w)))$. The dg category $\MF_{\hat{G}}(R,w)$ is pretriangulated, as can be verified directly. Since $\MF(R,w)$ is triangulated \cite[Lemma 5.6]{dyckerhoff2011}, we can apply \cite[Theorem 8.7]{elagin2014} to conclude that the canonical functor $H^0(\Perf(\MF_{\hat{G}}(R,w))) \rightarrow \HMF(R,w)^{\hat{G}}$ is an equivalence.
\end{proof}

\begin{remark}\label{rem:idemComp}
The definition of $\HMF_{G}(R,w)$ given in \cref{def:RealMF} agrees with \cite[\S 2.2]{ashok2004}, \cite[\S 6]{velez2009}, \cite[\S 3.1]{hori2008}, \cite[\S 2.1]{polishchuk2012}. On the other hand, Carqueville--Runkel, for example, define $\HMF_{G}(R,w)$ to be $\HMF(R,w)^G$ \cite[\S 7.1]{carqueville2016b}. In view of \cref{prop:idemCompHMFG}, these definitions agree for $w$ an isolated hypersurface singularity precisely if $\HMF_G(R,w)$ is idempotent complete, which does not seem to be known. If $\HMF_G(R,w)$ (or its Real generalization) is idempotent complete, then some of the proofs below can be simplified by using part \ref{prop:fixedPointFunctorQusiEquivOne} of \cref{prop:fixedPointFunctorQusiEquiv} in place of part \ref{prop:fixedPointFunctorQusiEquivTwo}.
    
For comparison, if $X$ is a quasi-projective $G$-variety over a field $k$ with $\ch \, k \nmid \vert G \vert$, then $D^b(\Coh(X)^G) \rightarrow D^b(\Coh(X))^G$ is a triangle equivalence \cite[Theorem 9.6]{elagin2011}, \cite[Theorem 7.3]{elagin2014}. A key ingredient of the proof is idempotent completeness of the bounded derived category of an abelian category \cite[Corollary 2.10]{balmer2001}.
\end{remark}
    
\subsection{Real Kn\"{o}rrer periodicity}
    
To formulate a Real generalization of Kn\"orrer periodicity (\cref{knoerreroriginal}), we study Real $G$-equivariant structures on $K=\{y+iz,y-iz\}$.
    
\begin{proposition} \label{prop:RealEquivRankOneMF}
The matrix factorization $\{u,v\} \in \MF(\mathbb{C}\pser{u,v},uv)$ admits a Real $G$-equivariant structure if and only if there exists $\chi \in Z^1(\hat{G};\mathbb{C}^{\times}_{\pi})$ such that $\sigma(u)=\chi(\sigma)u$ and $\sigma(v)=\chi(\sigma)^{-1} v$ for each $\sigma \in \hat{G}$, in which case the set of dg isomorphism classes of Real $G$-equivariant structures on $\{u,v\}$ is in bijection with $H^1(\hat{G}; \mathbb{C}\pser{u,v}^{\times}_{\pi})$.
\end{proposition}	
\begin{proof}
A Real $G$-equivariant structure on $\{u,v\}$ is the data of commutative diagrams
\begin{equation*}
\begin{tikzcd}
\mathbb{C}\pser{u,v} \arrow{r}{u} \arrow{d}[swap]{u_{\sigma}^{\chi_0}} & \mathbb{C}\pser{u,v} \arrow{d}[swap]{u_{\sigma}^{\chi_1}} \arrow{r}{v} & \mathbb{C}\pser{u,v} \arrow{d}{u_{\sigma}^{\chi_0}}\\
\mathbb{C}\pser{u,v}^{\sigma} \arrow{r}[swap]{u^{\sigma}} & \mathbb{C}\pser{u,v}^{\sigma} \arrow{r}[swap]{v^{\sigma}}& \mathbb{C}\pser{u,v},
\end{tikzcd}
\end{equation*}
where $\chi_i \in Z^1(\hat{G}; \mathbb{C}\pser{u,v}^{\times}_{\pi})$ and $u_{\sigma}^{\chi_i}(r)= \chi_i(\sigma^{-1}) \sigma^{-1}(r)$ for each $\sigma \in \hat{G}$ and $r \in \mathbb{C}\pser{u,v}$. Set $\chi=\chi_0 \chi_{1}^{-1}$. Commutativity of the diagram is then equivalent to the equations $\sigma(u)=\chi(\sigma)u$ and $\sigma(v)=\chi(\sigma)^{-1} v$. Since $\hat{G}$ acts on $\mathbb{C}\pser{u,v}$ by degree preserving maps, we find the restriction $\chi \in Z^1(\hat{G}; \mathbb{C}^{\times}_{\pi})$. Direct calculation then shows that isomorphism classes of Real $G$-equivariant structures are in bijection with
\[
\{(\chi_0,\chi_1) \in Z^1(\hat{G}; \mathbb{C}\pser{u,v}^{\times}_{\pi})^2 \mid \chi= \chi_0 \chi_1^{-1}\} \slash B^1(\hat{G}; \mathbb{C}\pser{u,v}^{\times}_{\pi}).
\]
Projection to the first factor gives the desired bijection with $H^1(\hat{G}; \mathbb{C}\pser{u,v}^{\times}_{\pi})$.
\end{proof}
		
By \cref{prop:RealEquivRankOneMF}, $\{u,v\}$ admits a Real $G$-equivariant structure precisely when $\mathbb{C}\pser{u,v} \simeq \widehat{\Sym\, V^{\vee}}$ with $V \simeq \mathbb{C}_{\chi^{-1}} \oplus \mathbb{C}_{\chi}$ for some $\chi \in Z^1(\hat{G};\mathbb{C}^{\times}_{\pi})$, where $\mathbb{C}_{\chi}$ is the one dimensional Real representation of $G$ determined by $\chi$ and $u$ (resp. $v$) is the coordinate dual to $\mathbb{C}_{\chi^{-1}}$ (resp. $\mathbb{C}_{\chi}$). Up to isomorphism, $\mathbb{C}_{\chi}$ depends on $\chi$ through $[\chi] \in H^1(\hat{G} ; \mathbb{C}^{\times}_{\pi})$. In the terminal case, $H^1(\Ctwo ; \mathbb{C}^{\times}_{\pi})$ is trivial so that, without loss of generality, we may take $\chi(\sigma)=-1$ for the generator $\sigma \in \Ctwo$. Pull back along the $\Ctwo$-grading then gives a universal choice of $\chi$ for all $\Ctwo$-graded finite groups $\hat{G}$. In what follows, fix a Real $G$-equivariant structure on $\{u,v\}$ by taking $\chi_0=\chi$ and $\chi_1$ to be trivial.

\cref{prop:RealEquivRankOneMF} applies to $\{y+iz , y-iz\} \in \MF(\mathbb{C}\pser{y,z},y^2+z^2)$ via the coordinate change $y = \frac{u+v}{2}$ and $z = \frac{u-v}{2i}$. We conclude that $\{y+iz , y-iz\}$ admits a Real structure precisely when, up to Real equivariant isomorphism, $\Ctwo$ acts on $\mathbb{C}\pser{y,z}$ by $(y,z) \mapsto (-y,z)$.
	
We can now state the first form of Real Kn\"{o}rrer periodicity.
	
\begin{theorem}\label{thm:RealKnorrer}
Let a $\Ctwo$-graded finite group $\hat{G}$ act on $R= \mathbb{C}\pser{x_1,\ldots, x_n}$ by generalized algebra automorphisms which preserve the potential $w$. Extend the generalized $\hat{G}$-action to $R\pser{y,z}$ by $\sigma(y)=\pi(\sigma) y$ and $\sigma(z)=z$, $\sigma \in \hat{G}$. Let $K= \{y+iz, y-iz\}$, considered as a Real $G$-equivariant matrix factorization as above. The dg functor
\begin{equation*}
\Ind \mathcal{K}^{\hat{G}} : \Perf(\MF_{\hat{G}}(R, w)) \rightarrow \Perf(\MF_{\hat{G}}(R\pser{y,z}, w + y^2+z^2))
\end{equation*}
is a quasi-equivalence.
\end{theorem}

\begin{proof}
By \cref{equivariantfunctor,lem:extTensorEquiv}, the Real $G$-equivariant structure on $K$ induces a Real $G$-equivariant structure on the Kn\"{o}rrer dg functor $\mathcal{K}$. 
The theorem therefore follows from \cref{knoerreroriginal} and \cref{prop:fixedPointFunctorQusiEquiv}\ref{prop:fixedPointFunctorQusiEquivOne}.
\end{proof}
    
\begin{corollary}
Let $w \in \mathbb{R}\pser{x_1, \dots, x_n}$ be an isolated hypersurface singularity at the origin. There is a triangle equivalence
$$ \HMF(\mathbb{R}\pser{x_1,\ldots,x_n},w) \simeq \HMF(\mathbb{R}\pser{x_1,\ldots,x_n,y,z}, w - y^2+z^2).$$
\end{corollary}
	
\begin{proof}
Let the terminal $\Ctwo$-graded group act on $\mathbb{C}\pser{x_1,\ldots,x_n}$ by complex conjugation and on $\mathbb{C}\pser{y,z}$ by the antilinear extension of $y \mapsto -y$ and $z \mapsto z$. The $\Ctwo$-fixed point subalgebra of $\mathbb{C}\pser{x_1,\ldots,x_n,y,z}$ is $\mathbb{R}\pser{x_1, \dots, x_n, iy,z}$, which is isomorphic to $\mathbb{R}\pser{x_1,\ldots,x_n,y,z}$. Under this isomorphism, the potential $w + y^2 +z^2$ is identified with $w-y^2+z^2$. In view of \cref{ex:antilinInvtPot}\ref{ex:antilinInvtPotOverR} and the fact that both triangulated categories under consideration are idempotent complete, \cref{thm:RealKnorrer} reduces to the desired statement.
\end{proof}
	
When the $\Ctwo$-grading of $\hat{G}$ is trivial, \cref{thm:RealKnorrer} gives a simple proof of the affine case of Hirano's equivariant Kn\"{o}rrer periodicity \cite[Theorem 1.2]{hirano2017}.
	
\begin{corollary}\label{cor:babyHirano}
Let a finite group $G$ act on $\mathbb{C}\pser{x_1, \dots, x_n}$ by $\mathbb{C}$-algebra automorphisms which preserve the potential $w$. There is a quasi-equivalence
\[
\Perf(\MF_G(\mathbb{C}\pser{x_1, \dots, x_n},w)) \simeq \Perf(\MF_G(\mathbb{C}\pser{x_1, \dots, x_n,y,z}, w + y^2+z^2)),
\]
where $G$ acts trivially on $\mathbb{C}\pser{y,z}$.
\end{corollary}
	
\begin{remark}\label{rem:charExtension}
There is a straightforward generalization of \cref{thm:RealKnorrer}, with the same proof, in which a Real character $\chi \in Z^1(\hat{G};\mathbb{C}^{\times}_{\pi})$ is used to extend the generalized $\hat{G}$-action from $R$ to $R \pser{y,z}$, as in \cref{prop:RealEquivRankOneMF}. When the $\Ctwo$-grading of $\hat{G}$ is trivial, this generalization recovers the extra character theoretic data in Hirano's equivariant Kn\"{o}rrer periodicity. 
\end{remark}

\section{Clifford modules and Real matrix factorizations} \label{sec:quadPot}

In this section, we study Real matrix factorizations of quadratic hypersurfaces using Clifford modules. In this way, we obtain Real generalizations of results of Buchweitz--Eisenbud--Herzog \cite{buchweitz1987} and connect Real Kn\"{o}rrer periodicity to classical periodicities of categories of Clifford modules \cite{atiyah1964}, \cite{lawson1989}.

\subsection{Buchweitz--Eisenbud--Herzog equivalence} \label{sec:BEH}

For background on Clifford algebras, see \cite[\S 1]{atiyah1964} and \cite[\S I]{lawson1989}; our conventions match the former.

Let $V$ be a finite dimensional vector space over a field $k$ with $\ch \, k \neq 2$. Let $R = \widehat{\Sym\, V^{\vee}}$ be the completed symmetric algebra on $V^{\vee}$ and $T(V)$ the tensor algebra on $V$. The Clifford algebra of a non-degenerate quadratic form $q \in \Sym^2\, V^{\vee}$ is the $\Ztwo$-graded algebra
\[
\Cl(V,q) = T(V) \slash ( v \otimes v - q(v)  \mid v \in V ).
\]
Let $\Cl(V,q) \mhyphen \grmod$ be the category of finitely generated $\Ztwo$-graded left $\Cl(V,q)$-modules. Objects are therefore $\Ztwo$-graded vector spaces $A=A_0 \oplus A_1$ with an action of $\Cl(V,q)$ such that homogeneous elements $c \in \Cl(V,q)$ act by linear maps $c: A_i \rightarrow A_{i + \vert c \vert}$. Morphisms are $\Cl(V,q)$-linear maps of degree zero.

Define a functor $\Phi: \Cl(V,q) \mhyphen \grmod \rightarrow Z^0(\MF(R,q))$ as follows \cite[\S 2]{buchweitz1987}. For $A \in \Cl(V,q) \mhyphen \grmod$, let $\Phi(A)=A \otimes_k R$ as $\Ztwo$-graded $R$-modules with twisted differential
\[
d_{\Phi(A)}^i \in \Hom_R(A_i \otimes_k R, A_{i+1} \otimes_k R) \simeq \Hom_k(A_i,A_{i+1}) \otimes_k R
\]
defined to be the image of the identity map under the composition
\[
\End_k(V) \xrightarrow[]{\sim} V \otimes_k V^{\vee} \rightarrow \Hom_k(A_i,A_{i+1}) \otimes_k R,
\]
where the second map uses the action of $V \subset \Cl(V,q)$ on $A$ and the canonical inclusion $V^{\vee} \hookrightarrow R$. Given a morphism $f: A \rightarrow A^{\prime}$ in $\Cl(V,q) \mhyphen \grmod$, set $\Phi(f) : \Phi(A) \xrightarrow[]{f \otimes \id_R} \Phi(A^{\prime})$.

\begin{theorem}[{\cite[\S 2]{buchweitz1987}}] \label{thm:BEH}
The functor $\Phi$ induces a $k$-linear equivalence
\[
\Cl(V,q) \mhyphen \grmod \xrightarrow[]{\sim} \HMF(\widehat{\Sym\, V^{\vee}},q).
\]
\end{theorem}

\subsection{Buchweitz--Eisenbud--Herzog equivalence and Reality}

Suppose now that $\hat{G}$ is a $\Ctwo$-graded finite group, $V$ is a finite dimensional Real representation of $G$ over $k=\mathbb{C}$ and the non-degenerate quadratic form $q \in \Sym^2 \, V^{\vee}$ is $\hat{G}$-invariant. Then $\hat{G}$ acts on $\Cl(V,q)$ by generalized algebra automorphisms so that, by \cref{2rep}, $\Cl(V,q) \mhyphen \grmod$ is a Real $2$-representation of $G$. Objects of $\Cl(V,q) \mhyphen \grmod^{\hat{G}}$ are called \emph{Real $G$-equivariant Clifford modules}. When $\hat{G} = \Ctwo$, this recovers the Real Clifford modules of \cite[\S 4]{atiyah1966}, \cite[\S I.10]{lawson1989}. By \cref{lem:Real2RepMF}, the dg category $\MF(R,q)$ is a Real $2$-representation of $G$, as is its subcategory $Z^0(\MF(R,w))$.

\begin{proposition} \label{prop:equivStrBEHFunctor}
The functor $\Phi: \Cl(V,q) \mhyphen \grmod \rightarrow Z^0(\MF(R,q))$ admits a Real $G$-equivariant structure.
\end{proposition}

\begin{proof}
Let $A \in \Cl(V,q) \mhyphen\grmod$ and $\sigma \in \hat{G}$. Define an $R$-linear isomorphism
\[
\eta_{\sigma,A} : \Phi(A^{\sigma}) \rightarrow \Phi(A)^{\sigma},
\qquad
a \otimes r
\mapsto
a \otimes \sigma^{-1}(r).
\]
That $\eta_{\sigma,A}$ is closed of degree zero is the statement that the diagrams
\begin{equation*}
		    \begin{tikzcd}[ampersand replacement=\&, column sep = 1cm, row sep = 1cm]
				A_i^{\sigma} \otimes R \arrow[r, "d^i"] \arrow[d, "\eta^i_{\sigma}", swap]\& A_{i+1}^{\sigma} \otimes R \arrow[d, "\eta^{i+1}_{\sigma}"]\\
				A_i \otimes R^{\sigma} \arrow[r, "d^{i \sigma}", swap] \& A_{i+1} \otimes R^{\sigma} 
			\end{tikzcd}
\end{equation*}
commute, $i=0,1$. The clockwise and counterclockwise compositions send $a \otimes r$ to $\sum_j \sigma^{-1} (v_j) a \otimes \sigma^{-1}(v^{\vee}_j r)$ and $\sum_j v^{\prime}_j a \otimes v_j^{\prime \vee} \sigma^{-1}(r)$, respectively, where $\{v_j\}$ and $\{v_j^{\prime}\}$ are arbitrary bases of $V$ with dual bases $\{v^{\vee}_j\}$ and $\{v_j^{\prime \vee}\}$. Taking $v_j^{\prime} = \sigma^{-1} (v_j)$
yields the required commutativity.
The equalities
\[
\eta^i_{\sigma_1} \circ \eta^i_{\sigma_2}(a \otimes r)
=
\eta^i_{\sigma_1}(a \otimes \sigma_2^{-1}(r))
=
a \otimes \sigma_1^{-1}(\sigma_2^{-1}(r))
=
\eta^i_{\sigma_2 \sigma_1}(a \otimes r)
\]
verify the coherence condition \eqref{dia:equivariantfunctor}.
\end{proof}

\begin{theorem}\label{realbeh}
Let $V$ be a Real representation of a finite $\Ctwo$-graded group $\hat{G}$ and $q \in \Sym^2\, V^{\vee}$ a non-degenerate $\hat{G}$-invariant quadratic form. There is an $\mathbb{R}$-linear equivalence
\[
\Cl(V,q)\mhyphen\grmod^{\hat{G}} \simeq \HMF(\widehat{\Sym\, V^{\vee}},q)^{\hat{G}}.
\]
\end{theorem}
	
\begin{proof}
Write $R = \widehat{\Sym\, V^{\vee}}$. The Real $2$-representation of $G$ on $\MF(R,w)$ induces one on $\HMF(R,w)$, for which the canonical functor $Z^0(\MF(R,w)) \rightarrow \HMF(R,w)$ admits an obvious Real $G$-equivariant structure. In view of \cref{prop:equivStrBEHFunctor}, the following composition
\[
\widetilde{\Phi}:\Cl(V,q) \mhyphen \grmod \xrightarrow[]{\Phi} Z^0(\MF(R,w)) \rightarrow \HMF(R,w)
\]
is Real $G$-equivariant. Moreover, $\widetilde{\Phi}$ is an equivalence by \cref{thm:BEH}. Applying \cite[Proposition 2.20]{sun2019}, we conclude that $\widetilde{\Phi}^{\hat{G}}:\Cl(V,q)\mhyphen\grmod^{\hat{G}} \rightarrow \HMF(R,q)^{\hat{G}}$ is an equivalence.
\end{proof}

\subsection{Periodicities of Clifford modules}\label{sec:cliffPeriod}
	
In this section, we connect \cref{thm:RealKnorrer,realbeh} to classical periodicities of categories of Clifford modules over $\mathbb{R}$.
	
Let $V=\bigoplus_{j=1}^n \mathbb{C} \cdot e_j$ with coordinates $x_j = e_j^{\vee}$ and $R = \widehat{\Sym \, V^{\vee}}$. Define a Real representation of the trivial group on $V$ by
\[
\sigma(e_j)=e_j, \;\; j=1, \dots, r, \qquad
\sigma(e_j)=-e_j, \;\; j=r+1, \dots, n
\]
where $\sigma \in \Ctwo$ is the generator. The quadratic form $q=\sum_{j=1}^n x_j^2$ is $\Ctwo$-invariant. Setting $s=n-r$, there is an equivalence $\Cl(V,q)\mhyphen\grmod^{\Ctwo} \simeq \Cl_{r,s} \mhyphen\grmod$, the right hand side being the category of graded modules over the $\mathbb{R}$-algebra $\Cl_{r,s} = \Cl(\mathbb{R}^n, \sum_{j=1}^r x_j^2 - \sum_{j=r+1}^n x_j^2)$. See \cite[Proposition I.10.9]{lawson1989}.
	
The categories $\Cl_{r,s} \mhyphen\grmod$, $r,s \in \mathbb{Z}_{\geq 0}$, enjoy a number of periodicities \cite[\S 4]{atiyah1964}, \cite[\S I.4]{lawson1989} which, in view of \cref{realbeh}, correspond to periodicities of Real matrix factorization categories. In particular, the $(1,1)$-periodicity theorem
\[
\Cl_{r,s} \mhyphen\grmod \simeq \Cl_{r+1,s+1} \mhyphen\grmod
\]
corresponds to an equivalence
\[
\HMF(R,q)^{\Ctwo} \simeq \HMF(R\pser{y,z},q+y^2+z^2)^{\Ctwo},
\]
where the generalized $\Ctwo$-action is extended from $R$ to $R\pser{y,z}$ by $\sigma(y)=-y$ and $\sigma(z)=z$. This is precisely the form of \cref{thm:RealKnorrer}, thereby making precise its $(1,1)$-nature. As a second example, recall that there are graded algebra isomorphisms
\[
\Cl_{1,1}^{\otimes 4}
\simeq
\Cl_{4,4}
\simeq
\Cl_{0,8}
\simeq
\Cl_{8,0}.
\]
In view of \cref{ex:antilinInvtPot}\ref{ex:antilinInvtPotOverR}, repeated application of \cref{thm:RealKnorrer}, combined with \cref{realbeh}, recovers a special case of Brown's $8$-periodicity.

\begin{corollary}[{\cite[Theorem 1.2]{brown2016}}] \label{cor:real8Periodicity}
Let $q$ be a non-degenerate quadratic form on $\mathbb{R}^n$. There is an $\mathbb{R}$-linear triangle equivalence
$$ \HMF(\mathbb{R}\pser{x_1,\ldots,x_n},q) \simeq \HMF(\mathbb{R}\pser{x_1,\ldots,x_{n+8}}, q + \sum_{i=1}^8x^2_{n+i}).$$
\end{corollary}

Note that Brown's theorem applies to arbitrary isolated hypersurface singularities, not only quadratics, but relies on Dyckerhoff's $2$-periodic derived Morita theory \cite[\S 6]{dyckerhoff2011}. On the other hand, our methods are elementary. It would be interesting to generalize Dyckerhoff's results to the Real equivariant setting and use them to deduce \cref{thm:RealKnorrer} from the Morita triviality of $\Cl_{1,1}$.
    
\section{Kn\"{o}rrer periodicity for Landau--Ginzburg orientifolds} \label{sec:orientifoldMF}

Motivated by Grothendieck--Witt theory and the physics of Landau--Ginzburg orientifolds, as studied by Hori--Walcher \cite{hori2008}, we give a second formulation of Real matrix factorizations and their periodicity, independent from those of \cref{sec:RealKnorrer,sec:quadPot}. All constructions are now linear over a ground field $k$ with $\ch\, k \neq 2$ and Real $2$-representations are in the contravariant approach of \cref{sec:contravariantApproach}.

\subsection{$2$-periodic dg categories with duality}
To begin this section, we define duality structures on ($2$-periodic) dg categories, parallel to the $\mathbb{Z}$-graded setting of \cite[\S 1]{schlichting2017}.

\begin{definition} \label{def:dgCatWD}
\begin{enumerate}[wide,labelwidth=!, labelindent=0pt,label=(\roman*)]
\item A \emph{dg category with duality} $(\catC, P, \Theta)$ is a dg category $\catC$, a dg functor $P: \catC^{\op} \rightarrow \catC$ and a dg natural isomorphism $\Theta: \id_{\catC} \Rightarrow P \circ P^{\op}$ such that
\begin{equation} \label{eq:doubleDualCoherence}
P(\Theta_C) \circ \Theta_{P(C)} = \id_{P(C)}
\end{equation}
for each $C \in \catC$. The pair $(P, \Theta)$ is called a \emph{dg duality structure} on $\catC$.
\item A \emph{dg form functor} $(T,\varphi): (\catC, P, \Theta) \rightarrow (\catD, Q, \Xi)$ between dg categories with duality is a dg functor $T: \catC \rightarrow \catD$ and a dg natural transformation $\varphi: T \circ P \Rightarrow Q \circ T^{\op}$ such that
\begin{equation}\label{eq:formFuncCoher}
Q(\varphi_C) \circ \Xi_{T(C)}
=
\varphi_{P(C)} \circ T(\Theta_C)
\end{equation}
for each $C \in \catC$.\qedhere
\end{enumerate}
\end{definition}

Given dg form functors $(T,\varphi): (\catC, P, \Theta) \rightarrow (\catD, Q, \Xi)$ and $(S,\psi): (\catD, Q, \Xi) \rightarrow (\catE, R, \Lambda)$, their composition is the dg form functor $(S \circ T, \psi \star \varphi): (\catC, P, \Theta) \rightarrow (\catE, R, \Lambda)$, where $\psi \star \varphi$ is defined so that its component at $C \in \catC$ is $\psi_{T(C)} \circ S(\varphi_C)$. Denote by $\dgCatD_k$ the category with objects dg categories with duality and morphisms their dg form functors.

\begin{example} \label{ex:Real2RepVsDuality}
Let $(\rho, \theta)$ be a Real $2$-representation of the trivial group (with respect to the terminal $\Ctwo$-graded group) on a dg category $\catC$. Let $\sigma$ be the generator of $\Ctwo$. Define a dg duality structure on $\catC$ by $P=\rho(\sigma)$ and $\Theta = \theta_{\sigma,\sigma}^{-1} \circ \theta_e^{-1}$. Similarly, a Real equivariant functor between Real $2$-representations induces a dg form functor. These assignments extend to an equivalence $\Ctwo \mhyphen \dgCat_k \xrightarrow[]{\sim} \dgCatD_k$.
\end{example}

Let $\hat{G}$ be a $\Ctwo$-graded finite group. We formulate two general results relating Real $2$-representations of $G$ and duality structures. In view of \cref{ex:Real2RepVsDuality}, the first is an instance of the expectation that, given an algebraic or geometric object with an action of a group $G$, the fixed points of a normal subgroup $N \unlhd G$ inherit an action of $G \slash N$. Our setting differs from standard ones (e.g. \cite[Exercise 4.15.3]{etingof2015}), since our actions are not covariant. Some calculations in the proof are relegated to Appendix \ref{app:A}.

\begin{theorem} \label{thm:dualityEquivCat}
Let $(\rho, \theta)$ be a Real $2$-representation of $G$ on a dg category $\catC$
\begin{enumerate}[wide,labelwidth=!, labelindent=0pt,label=(\roman*)]
\item Each $\sigma \in \hat{G} \setminus G$ induces a dg duality structure $(\rho(\sigma), \Theta(\sigma))$ on $\catC^G$.
\item \label{thm:dualityEquivCatUnique} For each $\sigma_1,\sigma_2 \in \hat{G} \setminus G$, there exist dg form equivalences
\[
(\id_{\catC^G},\varphi^{\sigma_1, \sigma_2}) : (\catC^G,\rho(\sigma_1), \Theta(\sigma_1)) \rightarrow (\catC^G,\rho(\sigma_2), \Theta(\sigma_2))
\]
which satisfy the coherence condition
\[
(\id_{\catC^G},\varphi^{\sigma_1 , \sigma_3}) = (\id_{\catC^G},\varphi^{\sigma_3, \sigma_2}) \circ (\id_{\catC^G},\varphi^{\sigma_1, \sigma_2}).
\]
\end{enumerate}
In particular, $\catC^G$ inherits from $(\rho,\theta)$ a canonical equivalence class of dg duality structures.
\end{theorem}

\begin{proof}
\begin{enumerate}[wide,labelwidth=!, labelindent=0pt,label=(\roman*)]
\item Lift $\rho(\sigma): \catC^{\op} \rightarrow \catC$ to a dg functor $(\catC^G)^{\op} \rightarrow \catC^G$, again denoted by $\rho(\sigma)$, as follows. Given $(C,u) \in \catC^G$, let $v=\{v_g\}_{g \in G}$, where
\[
v_g: \rho(\sigma)(C) \xrightarrow[]{\rho(\sigma)(u_{\sigma^{-1} g \sigma}^{-1})} \rho(\sigma) (\rho(\sigma^{-1} g \sigma)(C)) \xrightarrow[]{\theta_{\sigma, \sigma^{-1} g \sigma}} \rho(g \sigma)(C) \xrightarrow[]{\theta_{g,\sigma}^{-1}} (\rho(g) \rho(\sigma)(C)).
\]
\cref{lem:functorOnFixedPoints} verifies that $(\rho(\sigma)(C),v)$ is an object of $\catC^G$. Setting $\rho(\sigma)(C,u) = (\rho(\sigma)(C),v)$ and defining $\rho(\sigma)$ on morphisms in $\catC^G$ as for $\catC$ gives a dg functor $\rho(\sigma):(\catC^G)^{\op} \rightarrow \catC^G$. Let $\Theta(\sigma)$ be the natural isomorphism whose component at $(C,u) \in \catC^G$ is
\[
\Theta(\sigma)_{(C,u)}: C \xrightarrow[]{u_{\sigma^2}} \rho(\sigma^2)(C) \xrightarrow[]{\theta^{-1}_{\sigma,\sigma}} \rho(\sigma)^2(C).
\]
\cref{lem:ThetaMorphInCG} verifies that $\Theta(\sigma)_{(C,u)}$ is a morphism in $\catC^G$. We have
\[
\sigma(\Theta(\sigma)_{(C,u)}) \circ \Theta(\sigma)_{\sigma(C,u)}
=
\sigma(u_{\sigma^2}) \circ \sigma(\theta_{\sigma,\sigma}^{-1}) \circ \theta_{\sigma,\sigma}^{-1} \circ \theta_{\sigma^2, \sigma}^{-1} \circ \theta_{\sigma, \sigma^2} \circ \sigma(u^{-1}_{\sigma^2}),
\]
which is the identity by equation \eqref{eq:theta2CocycleContra}. Hence, $(\rho(\sigma),\Theta(\sigma))$ satisfies the coherence condition \eqref{eq:doubleDualCoherence} and defines a dg duality structure on $\catC^G$.

\item Let $\varphi^{\sigma_1,\sigma_2}: \rho(\sigma_1) \Rightarrow \rho(\sigma_2)$ be the natural isomorphism whose component at $(C,u) \in \catC^G$ is
\[
\varphi^{\sigma_1,\sigma_2}_{(C,u)}: \rho(\sigma_1)(C) \xrightarrow[]{\theta^{-1}_{\sigma_2, \sigma_2^{-1} \sigma_1}} \rho(\sigma_2)(\rho(\sigma_2^{-1} \sigma_1)(C)) \xrightarrow[]{\rho(\sigma_2)(u_{\sigma_2^{-1} \sigma_1})} \rho(\sigma_2)(C).
\]
\cref{lem:dgFormNatural} verifies that $\varphi^{\sigma_1, \sigma_2}_{(C,u)}$ is a morphism in $\catC^G$. By \cref{lem:dgFormCoher}, the pair $(\id_{\catC^G},\varphi^{\sigma_1,\sigma_2})$ is a dg form equivalence. The coherence condition involving $\sigma_1, \sigma_2,\sigma_3 \in \hat{G} \setminus G$ is proved in \cref{lem:dgFormFunctTorsor}.
\qedhere
\end{enumerate}
\end{proof}

\begin{corollary}\label{cor:dgDualityPerf}
A Real $2$-representation of $G$ on $\catC$ induces on $\Perf(\catC^G)$ a canonical equivalence class of dg duality structures.
\end{corollary}

\begin{proof}
By \cref{thm:dualityEquivCat}, the dg category $\catC^G$ inherits a dg duality structure, unique up to equivalence. Interpreting this dg duality structure as a Real $2$-representation of the trivial group on $\catC^G$, as in \cref{ex:Real2RepVsDuality}, we conclude from \cref{prop:Real2RepOnPerf} that $\Perf(\catC^G)$ inherits a dg duality structure which is unique up to equivalence.
\end{proof}

The second general result asserts naturality of the dg duality structures constructed in \cref{thm:dualityEquivCat}, and therefore also in \cref{cor:dgDualityPerf}.

\begin{theorem} \label{thm:equivFunctorContra}
Let $(F,\eta): \catC \rightarrow \catD$ be a Real $G$-equivariant dg functor. 
\begin{enumerate}[wide,labelwidth=!, labelindent=0pt,label=(\roman*)]
\item For each $\sigma \in \hat{G} \setminus G$, there is an induced dg form functor
\[
(F^G,\psi^{\sigma}): (\catC^G,\rho_{\catC}(\sigma),\Theta(\sigma)) \rightarrow (\catD^G,\rho_{\catD}(\sigma),\Theta(\sigma))
\]
whose second component $\psi^{\sigma}$ is a natural isomorphism.

\item For each $\sigma_1,\sigma_2 \in \hat{G} \setminus G$, there is a commutative diagram of dg form functors
\[
\begin{tikzcd}[column sep=5em,row sep=2.5em]
(\catC^G,\rho_{\catC}(\sigma_1),\Theta(\sigma_1)) \arrow{r}{(F^G,\psi^{\sigma_1})} \arrow{d}[left]{(\id_{\catC^G},\varphi^{\sigma_1,\sigma_2})} & (\catD^G,\rho_{\catD}(\sigma_1),\Theta(\sigma_1)) \arrow{d}{(\id_{\catD^G},\varphi^{\sigma_1,\sigma_2})}\\
(\catC^G,\rho_{\catC}(\sigma_2),\Theta(\sigma_2)) \arrow{r}[below]{(F^G,\psi^{\sigma_2})} & (\catD^G,\rho_{\catD}(\sigma_2),\Theta(\sigma_2)).
\end{tikzcd}
\]
\end{enumerate}
\end{theorem}

\begin{proof}
\begin{enumerate}[wide,labelwidth=!, labelindent=0pt,label=(\roman*)]
\item Define $\psi^{\sigma}: F^G \circ \rho_{\catC}(\sigma) \Rightarrow \rho_{\catD}(\sigma) \circ (F^G)^{\op}$ so that its component at $(C,u) \in \catC^G$ is $\eta_{\sigma}: F(\rho_{\catC}(\sigma)C) \rightarrow \rho_{\catD}(\sigma) F(C)$. That $\psi^{\sigma}_{(C,u)}$ is a morphism in $\catD^G$ amounts to commutativity of the square
\begin{equation*}
\begin{tikzcd}[column sep=large,row sep=normal]
F(\rho_{\catC}(\sigma) C) \arrow[r, "\eta_{\sigma}"] \arrow[swap, d, "\eta_g \circ F(v^{\sigma}_g)"] & \rho_{\catD}(\sigma) F(C) \arrow[d, "\theta_{g,\sigma}^{-1} \circ \theta_{\sigma, \sigma^{-1} g \sigma} \circ \rho_{\catD}(\sigma)(F(u^{-1}_{\sigma^{-1} g \sigma}) \circ \eta^{-1}_{\sigma^{-1} g \sigma})"]\\
\rho_{\catD}(g) F(\rho_{\catC}(\sigma) C) \arrow[r, swap, "\rho_{\catD}(\eta_{\sigma})"] & \rho_{\catC}(g)(\rho_{\catD}(\sigma) F(C)))
\end{tikzcd}
\end{equation*}
for each $g \in G$.
The clockwise composition is
\begin{eqnarray*}
&&
\theta_{g,\sigma}^{-1} \circ \theta_{\sigma, \sigma^{-1} g \sigma} \circ \rho_{\catD}(\sigma)(\eta^{-1}_{\sigma^{-1} g \sigma}) \circ \rho_{\catD}(\sigma)(F(u^{-1}_{\sigma^{-1} g \sigma})) \circ \eta_{\sigma}\\
&=&
\theta_{g,\sigma}^{-1} \circ \theta_{\sigma, \sigma^{-1} g \sigma} \circ \rho_{\catD}(\sigma)(\eta^{-1}_{\sigma^{-1} g \sigma}) \circ \eta_{\sigma} \circ F(\rho_{\catC}(\sigma)(u^{-1}_{\sigma^{-1}g\sigma}))\\
&=&
\theta_{g,\sigma}^{-1} \circ \eta_{g \sigma} \circ F(\theta_{\sigma, \sigma^{-1} g \sigma}) \circ F(\rho_{\catC}(\sigma)(u^{-1}_{\sigma^{-1}g\sigma}))\\
&=&
\theta_{g, \sigma}^{-1} \circ \eta_{g \sigma} \circ F(\theta_{g, \sigma}) \circ F(\theta^{-1}_{g, \sigma}) \circ F(\theta_{\sigma, \sigma^{-1} g \sigma}) \circ F(\rho_{\catC}(\sigma)(u^{-1}_{\sigma^{-1}g \sigma})).
\end{eqnarray*}
The first equality follows from naturality of $\eta_{\sigma}$ and the second and third from equation \eqref{eq:equivariantFunctorContra}. The final expression is the counterclockwise composition.

To check that $(F^G, \psi^{\sigma})$ is a dg form functor, we verify equation \eqref{eq:formFuncCoher}, which reads
\[
\rho_{\catD}(\sigma)(\psi^{\sigma}_{(C,u)}) \circ \Theta(\sigma)_{F^G(C,u)}
=
\psi^{\sigma}_{\rho_{\catC}(\sigma)(C,u)} \circ F^G(\Theta(\sigma)_{(C,u)})
\]
for each $(C,u) \in \catC^G$. We compute
\begin{eqnarray*}
\rho_{\catD}(\sigma)(\psi^{\sigma}_{(C,u)}) \circ \Theta(\sigma)_{F^G(C,u)}
&=&
\rho_{\catD}(\sigma)(\eta_{\sigma}) \circ \theta^{-1}_{\sigma,\sigma} \circ \eta_{\sigma^2} \circ F(u_{\sigma^2}) \\
&=&
\eta_{\sigma} \circ F(\theta^{-1}_{\sigma,\sigma}) \circ F(u_{\sigma^2}) \\
&=&
\psi^{\sigma}_{\rho_{\catC}(\sigma)(C,u)} \circ F^G(\Theta(\sigma)_{(C,u)}),
\end{eqnarray*}
the second equality following from equation \eqref{eq:equivariantFunctorContra} and the others by definition.

\item The component at $(C,u) \in \catC^G$ of the counterclockwise composition is
\[
\eta_{\sigma_2} \circ F(\rho_{\catC}(\sigma_2)(u^{-1}_{\sigma_2^{-1} \sigma_1})) \circ F(\theta^{-1}_{\sigma_2,\sigma_2^{-1} \sigma_1})
=
\rho_{\catD}(\sigma_2)(F(u^{-1}_{\sigma_2^{-1} \sigma_1})) \circ \eta_{\sigma_2} \circ F(\theta^{-1}_{\sigma_2,\sigma_2^{-1} \sigma_1}),
\]
the equality following by naturality of $\eta_{\sigma_2}$. Using the definition of $F(C,u)$, the previous expression is seen to equal
\[
\rho_{\catD}(\sigma_2)(\eta_{\sigma_2^{-1} \sigma_1} \circ F(u^{-1}_{\sigma_2^{-1} \sigma_1})) \circ \theta^{-1}_{\sigma_2,\sigma_2^{-1} \sigma_1} \circ \eta_{\sigma_1}.
\]
Contravariance of $\rho_{\catD}(\sigma_2)$ and equation \eqref{eq:equivariantFunctorContra} with $\sigma_1 \mapsto \sigma_2^{-1} \sigma_1$ and $\sigma_2 \mapsto \sigma_2$ then gives equality with the clockwise composition.
\qedhere
\end{enumerate}
\end{proof}

\subsection{Duality for matrix factorizations} \label{sec:mfDuality}

In this section, we record standard material about duality for matrix factorizations, following \cite{yoshino1998}, \cite{polishchuk2012}.

Fix a potential $w \in R=k\pser{x_1,\dots,x_n}$. Given $M \in \MF(R,w)$, put $M_- = (M,-d_M) \in \MF(R,w)$. The grading morphism $J_M = \id_{M_0} \oplus (-\id_{M_1}) :  M \rightarrow M_-$ is a dg isomorphism which satisfies $J_{\Sigma M} = - \Sigma J_M$.

Let $(-)^{\vee} : \MF(R,w)^{\op} \rightarrow \MF(R,-w)$ be the dg functor defined on objects by $(M,d_M)^{\vee} = (M^{\vee}, d_{M^{\vee}})$, where $M^{\vee}=\Hom_R(M,R)$ and $d_{M^{\vee}}(f) = -(-1)^{\vert f \vert} f \circ d_M$, and on morphisms by $\Hom_R(-,R)$. In diagrammatic form, we have
\[
(\xymatrix@C+0.5pc{
M_0 \ar@<+.5ex>[r]^(0.5){d_M^0} & \ar@<+.5ex>[l]^(0.5){d_M^1}
M_1} )^{\vee} =  (\xymatrix@C+0.5pc{
M_0^{\vee} \ar@<+.5ex>[r]^(0.5){-d_M^{1 \vee}} & \ar@<+.5ex>[l]^(0.5){d_M^{0 \vee}} M_1^{\vee}}).
\]
The dg isomorphism $\Theta_M= J_{M^{\vee \vee}} \circ \ev_M : M \rightarrow M^{\vee \vee}$, where $\ev_M$ is the canonical evaluation $R$-module isomorphism $M \xrightarrow[]{\sim} M^{\vee \vee}$, makes $(\MF(R,w), (-)^{\vee}, \Theta)$ a dg category with duality. The grading morphisms are the components of a natural isomorphism
\begin{equation}\label{eq:shiftWithDual}
J: \Sigma \circ (-)^{\vee} \Rightarrow (-)^{\vee} \circ \Sigma.
\end{equation}

Let $M \in \MF(R,w)$ and $N \in \MF(R^{\prime},w^{\prime})$. Under the swap identification $R \otimes_k R^{\prime} \simeq R^{\prime} \otimes_k R$, there is an isomorphism
\begin{equation}\label{eq:tensorDual}
(M \boxtimes N)^{\vee} \simeq N^{\vee} \boxtimes M^{\vee}
\end{equation}
in $\MF(R \otimes_k R^{\prime}, -w \otimes 1 -1 \otimes w^{\prime})$ given by the pairing $\langle n^{\vee} \otimes m^{\vee}, m \otimes n \rangle=n^{\vee}(n) m^{\vee}(m)$.

\subsection{Orientifold data for $G$-equivariant matrix factorizations}

Let $\hat{G}$ be a $\Ctwo$-graded finite group acting on $R$ by $k$-algebra automorphisms. We stress that $\hat{G}$ acts by $k$-linear automorphisms, unlike the $\hat{G}$-actions of  \cref{sec:RealKnorrer,sec:quadPot}. Assume that the potential $w$ is $\pi$-semi-invariant:
\begin{equation*}
\sigma(w) = \pi(\sigma) w,
\qquad
\sigma \in \hat{G}.
\end{equation*}
The ungraded $G$ is therefore a group of symmetries of $w$, while $\hat{G}$ itself is not. The results of Section \ref{sec:RealMFCat}, applied to trivially $\Ctwo$-graded groups, imply that each $\sigma \in \hat{G}$ defines a dg functor $(-)^{\sigma} : \MF(R,w) \rightarrow \MF(R, \pi(\sigma)w)$. The identity maps are the components of a natural isomorphism
\begin{equation}\label{eq:dualGroupCommute}
(-)^{\sigma} \circ (-)^{\vee} \xRightarrow[]{\sim} (-)^{\vee} \circ ((-)^{\sigma})^{\op} : \MF(R,w)^{\op} \rightarrow \MF(R,\pi(\sigma)w).
\end{equation}
Given $g \in G$ and $\sigma \in \hat{G} \setminus G$, define dg functors $\rho(g): \MF(R,w) \xrightarrow[]{(-)^g} \MF(R,w)$ and
\[
\rho(\sigma): \MF(R,w)^{\op} \xrightarrow[]{(-)^{\vee}} \MF(R,-w) \xrightarrow[]{(-)^{\sigma}} \MF(R,w).
\]

\begin{lemma} \label{lem:contraReal2RepMF}
The dg functors $\{\rho(\sigma)\}_{\sigma \in \hat{G}}$ extend to a Real $2$-representation on $\MF(R,w)$.
\end{lemma}

\begin{proof}
Define coherence data $\{\theta_{\sigma_2,\sigma_1}\}_{\sigma_1, \sigma_2 \in \hat{G}}$ as follows. If $\pi(\sigma_2)=-1$ and $\pi(\sigma_1)=1$, let
\begin{multline*}
\theta_{\sigma_2,\sigma_1}: \rho(\sigma_2) \circ {^{\pi(\sigma_2)}}\rho(\sigma_1)
=
(-)^{\sigma_2} \circ (-)^{\vee} \circ ((-)^{\sigma_1})^{\op}
\xRightarrow[]{\sim} \\
(-)^{\sigma_2} \circ (-)^{\sigma_1} \circ (-)^{\vee} 
\xRightarrow[]{\sim}
(-)^{\sigma_2 \sigma_1}  \circ (-)^{\vee}
=
\rho(\sigma_2 \sigma_1),
\end{multline*}
where the first and second natural isomorphisms are those of \cref{eq:dualGroupCommute} and the proof of \cref{2rep} (restricted to ungraded groups), respectively. If $\pi(\sigma_2)= \pi(\sigma_1)=-1$, let
\begin{multline*}
\theta_{\sigma_2,\sigma_1}: 
\rho(\sigma_2) \circ {^{\pi(\sigma_2)}}\rho\sigma_1)
=
(-)^{\sigma_2} \circ (-)^{\vee} \circ ((-)^{\sigma_1} \circ (-)^{\vee})^{\op}
\xRightarrow[]{\sim} \\
(-)^{\sigma_2} \circ (-)^{\sigma_1} \circ (-)^{\vee} \circ ((-)^{\vee})^{\op}
\xRightarrow[]{\sim}
(-)^{\sigma_2 \sigma_1}
=
\rho(\sigma_2 \sigma_1),
\end{multline*}
where now the second natural isomorphism uses in addition the double dual isomorphism $\Theta$ for $\MF(R,w)$. The remaining $\theta_{\sigma_2,\sigma_1}$ are defined similarly. The coherence of $\theta$ follows from the coherence of \cref{2rep} and the coherence isomorphisms \eqref{eq:doubleDualCoherence} for $\Theta$.
\end{proof}

\begin{theorem}\label{thm:dualityEquivMF}
The Real $2$-representation of $G$ on $\MF(R,w)$ from \cref{lem:contraReal2RepMF} induces a canonical equivalence classes of dg duality structures on $\MF_G(R,w)$.
\end{theorem}

\begin{proof}
The follows immediately from \cref{thm:dualityEquivCat}.
\end{proof}

\begin{example}
Assume that the ground field $k$ is algebraically closed. For $m \in \mathbb{Z}_{\geq 1}$, denote by $\zeta_m \in k^{\times}$ a primitive $m\textsuperscript{th}$ root of unity and $C_m$ the cyclic multiplicative group of order $m$.
\begin{enumerate}[wide,labelwidth=!, labelindent=0pt,label=(\roman*)]
\item Let $w = x^m \in k\pser{x}$, $m \geq 2$, and $\hat{G} = C_{2m}$ with $\pi$ the projection with kernel $G= C_m$. Then $w$ is $\pi$-semi-invariant with respect to the $\hat{G}$-action on $k\pser{x}$ in which a generator acts by $x \mapsto \zeta_{2m} x$.

\item Let $w = x y^{2m} + x^{2n+1} \in k \pser{x,y}$, $m,n \in \mathbb{Z}_{\geq 1}$, and $\hat{G} = C_{2m}$ with $\pi$ the projection with kernel $G = C_m$. Then $w$ is $\pi$-semi-invariant with respect to the $\hat{G}$-action on $k\pser{x,y}$ in which a generator acts by $x \mapsto -x$ and $y \mapsto \zeta_{2m} y$.

\item Let $w = x^m - y^m \in k \pser{x,y}$, $m \in \mathbb{Z}_{\geq 2}$. As a first example, let $\hat{G} = C_m \times \Ctwo$ with $\pi$ the projection to the second factor. Define a $\hat{G}$-action on $k \pser{x,y}$ by letting a generator of $G$ act by $x \mapsto \zeta_m x$ and $y \mapsto \zeta_m y$ and a generator of $\Ctwo$ act by swapping $x$ and $y$. Then $w$ is $\pi$-semi-invariant. As a second example, let $\hat{G} = D_{2m}$ be the dihedral group of order $2m$ with $\pi$ the projection with kernel $G=C_m$. Define a $\hat{G}$-action on $k \pser{x,y}$ by letting a generator of $G$ act by $x \mapsto \zeta_m x$ and $y \mapsto \zeta_m^{-1} y$ and a fixed element of the non-identity coset of $\hat{G}$ act by swapping $x$ and $y$. Then $w$ is $\pi$-semi-invariant. Finally, let a generator of $\hat{G} = C_{2m}$ act on $k \pser{x,y}$ by $x \mapsto \zeta_{2m} x$ and $y \mapsto \zeta_{2m} y$. Then $w$ is again $\pi$-semi-invariant. The final example generalizes to a homogeneous potential in any number of variables in an obvious way. \qedhere
\end{enumerate}
\end{example}

There is a shifted variant of \cref{lem:contraReal2RepMF}, and so also \cref{thm:dualityEquivMF}, in which $\hat{G}$ acts on $\MF(R,w)$ by the dg functors $\tilde{\rho}(g) = \rho(g)$, $g \in G$, and $\tilde{\rho}(\sigma) = \rho(\sigma) \circ \Sigma$, $\sigma \in \hat{G} \setminus G$. The shifted coherence data $\tilde{\theta}$ is defined similarly to $\theta$, using in addition that $\Sigma$ commutes with $(-)^{\sigma}$, the natural isomorphism \eqref{eq:shiftWithDual} and $\Sigma^2 \simeq \id_{\MF}$.
For example, when $\pi(\sigma_2)= \pi(\sigma_1)=-1$, let
\begin{multline*}
\tilde{\theta}_{\sigma_2,\sigma_1}:
\tilde{\rho}(\sigma_2) \circ {^{\pi(\sigma_2)}}\tilde{\rho}(\sigma_1)
=
(-)^{\sigma_2} \circ (-)^{\vee} \circ \Sigma \circ ((-)^{\sigma_1} \circ (-)^{\vee} \circ \Sigma)^{\op}
\xRightarrow[]{\sim} \\
(-)^{\sigma_2} \circ (-)^{\sigma_1} \circ (-)^{\vee} \circ ((-)^{\vee})^{\op} \circ \Sigma \circ \Sigma
\xRightarrow[]{\sim}
(-)^{\sigma_2 \sigma_1}
=
\tilde{\rho}(\sigma_2 \sigma_1),
\end{multline*}
where now the grading isomorphism $J$ is used in the first natural isomorphism.

\begin{example}\label{ex:LGModelDiscTor}
As discussed in Section \ref{sec:twistCoherenceData}, a $2$-cocycle $\hat{\mu} \in Z^2(\hat{G}; k^{\times}_{\pi})$ can be used to twist the Real $2$-representation $(\rho,\theta)$, or $(\tilde{\rho},\tilde{\theta})$, of $G$ on $\MF(R,w)$. This construction admits the following physical interpretation. Let $\mu \in Z^2(G;k^{\times})$ be the underlying untwisted $2$-cocycle. As explained in \cite[\S 2.3]{brunner2015}, the homotopy fixed point category $\MF_{G,\mu}(R,w)$ of the $\mu$-twisted $2$-representation of $G$ consists of $D$-branes in the Landau--Ginzburg $G$-orbifold of $(R,w)$ with discrete torsion $\mu$. By \cref{thm:dualityEquivCat}, the data $(\hat{G}, \hat{\mu})$ defines on $\MF_{G,\mu}(R,w)$ a dg duality structure which, by generalizing the discussion in \cite[\S 3]{hori2008}, is data required to orientifold the above Landau--Ginzburg $G$-orbifold. In particular, objects of the $\Ctwo$-homotopy fixed point category $\MF_{G,\mu}(R,w)^{C_2} \simeq \MF_{\hat{G},\hat{\mu}}(R,w)$ are $D$-branes in the Landau--Ginzburg $\hat{G}$-orientifold of $(R,w)$ with discrete torsion $\hat{\mu}$.
\end{example}

\begin{remark}
An analogue of \cref{realbeh} holds in the contravariant setting. Since it is not used in what follows, we restrict ourselves to a brief discussion of the terminal case $\hat{G}=\Ctwo$. In the present setting, the generator $\sigma \in \Ctwo$ is required to negate the non-degenerate quadratic form, $\sigma(q)=-q$, and so yields a graded algebra involution
\[
\tilde{\sigma}: \Cl(V,q)^{\op} \xrightarrow[]{\sim} \Cl(V,-q) \xrightarrow[]{\sigma} \Cl(V,q),
\]
where $\Cl(V,q)^{\op}$ denotes the superalgebra opposite of $\Cl(V,q)$. The $k$-linear dual of a graded $\Cl(V,q)$-module is naturally a right $\Cl(V,q)$-module, which we view as a left $\Cl(V,q)$-module via $\tilde{\sigma}$. In this way, we obtain a duality structure on $\Cl(V,q) \mhyphen \grmod$ with respect to which we verify that $\tilde{\Phi}$ lifts to a form equivalence.
\end{remark}

\subsection{Real Kn\"{o}rrer periodicity}
    
Let $k$ be algebraically closed with $\ch\, k \nmid 2$ and $\ch \, k \nmid \vert G \vert$. We begin this section with a contravariant analogue of \cref{prop:RealEquivRankOneMF}.

\begin{proposition}\label{prop:RealEquivRankOneMFContra}
The matrix factorization $\{u,v\} \in \MF(k\pser{u,v},uv)$ admits a Real $G$-equivariant structure	\begin{enumerate}[wide,labelwidth=!, labelindent=0pt,label=(\roman*)]
\item \label{prop:RealEquivRankOneMFContraShifted} with respect to the Real $2$-representation $(\tilde{\rho},\tilde{\theta})$ if and only if there exists $\chi \in Z^1(\hat{G};k^{\times}_{\pi})$ such that $\sigma(u)=\pi(\sigma) \chi(\sigma) u$ and $\sigma(v)=\chi(\sigma)^{-1} v$ for each $\sigma \in \hat{G}$, and
	
\item \label{prop:RealEquivRankOneMFContraUnshifted} with respect to the Real $2$-representation $(\rho,\theta)$ if and only if there exists $\chi \in Z^1(\hat{G};k^{\times}_{\pi})$ such that
\[
\sigma \left( \begin{matrix} u \\ v \end{matrix} \right)
=
\begin{cases}
\left( \begin{smallmatrix} \chi(\sigma)u \\ \chi(\sigma)^{-1}v \end{smallmatrix} \right) & \mbox{if } \pi(\sigma)=1,\\
\left( \begin{smallmatrix} -\chi(\sigma)v \\ \chi(\sigma)^{-1}u \end{smallmatrix} \right) & \mbox{if } \pi(\sigma)=-1.
\end{cases}
\]
\end{enumerate}
In both cases, if non-empty, the set of dg isomorphism classes of Real $G$-equivariant structures on $\{u,v\}$ is in bijection with $H^1(\hat{G}; k\pser{u,v}^{\times}_{\pi})$.
\end{proposition}
	
\begin{proof}
\begin{enumerate}[wide,labelwidth=!, labelindent=0pt,label=(\roman*)]
\item A Real $G$-equivariant structure on $\{u,v\}$ is the data of commutative diagrams
\[
\begin{tikzcd}
R \arrow{r}{u} \arrow{d}[swap]{u_g^{\chi_0}} & R \arrow{d}[swap]{u_g^{\chi_1}} \arrow{r}{v} & R \arrow{d}{u_g^{\chi_0}}\\
R^g \arrow{r}[swap]{u^g} & R^g \arrow{r}[swap]{v^g}& R^g
\end{tikzcd}
\]
and
\[
\begin{tikzcd}
R \arrow{r}{u} \arrow{d}[swap]{u_{\sigma}^{\chi_0}} & R \arrow{d}[swap]{u_{\sigma}^{\chi_1}} \arrow{r}{v} & R \arrow{d}{u_{\sigma}^{\chi_0}}\\
R^{\vee \sigma} \arrow{r}[swap]{u^{\vee \sigma}} & R^{\vee \sigma} \arrow{r}[swap]{-v^{\vee \sigma}}& R^{\vee \sigma}
\end{tikzcd}
\]
where $R= k \pser{u,v}$, $g \in G$ and $\sigma \in \hat{G} \setminus G$. Here $\chi_i \in Z^1(\hat{G}; k\pser{u,v}^{\times}_{\pi})$. Setting $\chi=\chi_0 \chi_{1}^{-1}$, commutativity of the diagrams becomes the stated conditions and we conclude that $\chi \in Z^1(\hat{G}; k^{\times}_{\pi})$. The remainder of the proof is as for \cref{prop:RealEquivRankOneMF}.
	
\item The proof is a minor variation of the previous part and so is ommited. \qedhere
\end{enumerate}
\end{proof}
	
\cref{prop:RealEquivRankOneMFContra} distinguishes the Real $2$-representation $(\tilde{\rho}, \tilde{\theta})$: for any $\Ctwo$-graded finite group $\hat{G}$, there exists a Real $G$-equivariant structure on $\{u,v\}$ for some action of $\hat{G}$, say, by taking $\chi$ trivial. This is in contrast to $(\rho, \theta)$, where there exists no Real structure in the terminal case $\hat{G}=\Ctwo$. To formulate a universal form of Real Kn\"{o}rrer periodicity, we therefore work in the setting of \cref{prop:RealEquivRankOneMFContra}\ref{prop:RealEquivRankOneMFContraShifted} and take $k\pser{u,v} \simeq \widehat{\Sym\, V^{\vee}}$, where $V \simeq k_{\chi^{-1}} \oplus k_{\chi}$ for some $\chi \in Z^1(\hat{G};k^{\times}_{\pi})$ with $u$ (resp. $v$) the coordinate dual to $k_{\chi^{-1}}$ (resp. $k_{\chi}$). In the terminal case, the group $H^1(\Ctwo ; k^{\times}_{\pi})$ is trivial so that we may assume $\chi$ to be trivial. This gives a universal choice of $\chi$ for all $\Ctwo$-graded finite groups $\hat{G}$. Fix a Real $G$-equivariant structure on $\{u,v\}$, say, by taking $\chi_0=\chi_1$ trivial.
	
The next result uses the universal $2$-cocycle twist of a Real $2$-representation $(\rho,\theta)$, denoted by $(\rho, \theta_-)$ in Section \ref{sec:twistCoherenceData}, and views Kn\"{o}rrer periodicity (\cref{knoerreroriginal}) as a quasi-equivalence $\mathcal{K}= - \boxtimes \{u,v\} : \MF(R, w) \rightarrow \MF(R\pser{u,v}, w + uv)$.
	
\begin{theorem} \label{thm:oriKnorrer}
Let a $\Ctwo$-graded finite group $\hat{G}$ act on $R= \mathbb{C}\pser{x_1,\ldots, x_n}$ by $\mathbb{C}$-algebra automorphisms which leave the potential $w$ $\pi$-semi-invariant. Extend the $\hat{G}$-action to $R\pser{u,v}$ by $\sigma(u)=\pi(\sigma)u$ and $\sigma(v)=v$, $\sigma \in \hat{G}$. The Kn\"{o}rrer dg functor $\mathcal{K}$ admits a Real $G$-equivariant structure when $\MF(R,w)$ and $\MF(R\pser{u,v},w + uv)$ are viewed as Real $2$-representations by $(\rho, \theta)$ and $(\tilde{\rho}, \tilde{\theta}_-)$, respectively.
\end{theorem}
 
\begin{proof}
Define natural isomorphisms $\{\eta_{\sigma} : \mathcal{K} \circ \rho(\sigma) \Rightarrow \tilde{\rho}(\sigma) \circ {^{\pi(\sigma)}}\mathcal{K}\}_{\sigma \in \hat{G}}$ as follows. Let $M \in \MF(R,w)$. For $g \in G$, define the component of $\eta_g$ at $M$ to be
\[
\mathcal{K}(\rho(g)(M))
=
\rho(g)(M) \boxtimes K
\xrightarrow[]{\id \boxtimes u_g}
\rho(g)(M) \boxtimes \rho(g)(K)
\xrightarrow[]{\sim}
\rho(g)(M \boxtimes K),
\]
where the final map is the monoidal coherence data for $(-)^g$. For $\sigma \in \hat{G} \setminus G$, define the component of $\eta_{\sigma}$ at $M$ to be
\begin{multline*}
\mathcal{K}(\rho(\sigma)(M))
=
\rho(\sigma)(M) \boxtimes K
\xrightarrow[]{\sim}
K \boxtimes \rho(\sigma)(M)
\xrightarrow[]{u_{\sigma} \boxtimes \id}
\tilde{\rho}(\sigma)(K) \boxtimes \rho(\sigma)(M)\\
\xrightarrow[]{\sim}
\tilde{\rho}(\sigma)(M \boxtimes K) = \tilde{\rho}(\sigma)(\mathcal{K}(M)).
\end{multline*}
The first and third isomorphisms use the symmetry \eqref{eq:tensorSym} and the isomorphisms \eqref{eq:shiftTensorCompat} and \eqref{eq:tensorDual}, respectively. Explicitly, $\eta_{\sigma}$ is the morphism
\begin{equation}\label{eq:oddCompRealKnorrerContra}
\xymatrixcolsep{11pc}
\xymatrixrowsep{4pc}
\xymatrix{
M_0^{\vee \sigma} \oplus M_1^{\vee \sigma} \ar@<+.5ex>[r]^(0.5){\left( \begin{smallmatrix} -(d_M^1)^{\vee \sigma} & -v  \\ u & (d_M^0)^{\vee \sigma} \end{smallmatrix}\right)} \ar@<+.5ex>[d]_(0.5){\left( \begin{smallmatrix} 0 & 1 \\ -1 & 0 \end{smallmatrix}\right)} & \ar@<+.5ex>[l]^(0.5){\left( \begin{smallmatrix} (d_M^0)^{\vee \sigma} & v \\ -u & -(d_M^1)^{\vee \sigma} \end{smallmatrix}\right)}
M_1^{\vee \sigma} \oplus M_0^{\vee \sigma} \ar@<+.5ex>[d]^(0.5){\left( \begin{smallmatrix} 0 & 1 \\ 1 & 0 \end{smallmatrix}\right)}\\
M_1^{\vee \sigma} \oplus M_0^{\vee \sigma} \ar@<+.5ex>[r]^(0.5){\left( \begin{smallmatrix} (d_M^0)^{\vee \sigma} & -u  \\ -v & (d_M^1)^{\vee \sigma} \end{smallmatrix}\right)} & \ar@<+.5ex>[l]^(0.5){-\left( \begin{smallmatrix} (d_M^1)^{\vee \sigma} & u \\ v & (d_M^0)^{\vee \sigma} \end{smallmatrix}\right)}
M_0^{\vee \sigma} \oplus M_1^{\vee \sigma} .}
\end{equation}

It is immediate that each $\eta_{\sigma}$ is a natural isomorphism. Given the explicit description of $\eta_{\sigma}$ (for example, in the above matrix form), verification of the coherence condition \eqref{eq:equivariantFunctorContra}, with $\theta_{\catC} \mapsto \theta$ and $\theta_{\catD} \mapsto \tilde{\theta}_-$, is a direct calculation. For example, when $\pi(\sigma_2)=\pi(\sigma_1)=-1$ the coherence condition requires the equality
\[
	\tilde{\theta}_{-,\sigma_2,\sigma_1} \circ \left( \id_{\tilde{\rho}(\sigma_2)} \circ \eta_{\sigma_1}^{-1 \vee} \right) \circ \left(\eta_{\sigma_2} \circ \id_{{^{\pi(\sigma_2)}}\rho(\sigma_1)}\right)
	=
	\eta_{\sigma_2 \sigma_1} \circ \mathcal{K}(\theta_{\sigma_2, \sigma_1}).
\]
The components of the left hand side of the desired equality in degrees $0$ and $1$ are represented by the matrices $- \left(\begin{smallmatrix} 0 & 1 \\ -1 & 0 \end{smallmatrix} \right)^{- \textnormal{tr}} \left(\begin{smallmatrix} 0 & 1 \\ -1 & 0 \end{smallmatrix} \right) = \left(\begin{smallmatrix} 1 & 0 \\ 0 & 1 \end{smallmatrix}\right)$ and $(-1)^2 \left(\begin{smallmatrix} 0 & 1 \\ 1 & 0 \end{smallmatrix} \right)^{- \textnormal{tr}} \left(\begin{smallmatrix} 0 & 1 \\ 1 & 0 \end{smallmatrix} \right) = \left(\begin{smallmatrix} 1 & 0 \\ 0 & 1 \end{smallmatrix}\right)$, respectively, where we work in a basis as in diagram \eqref{eq:oddCompRealKnorrerContra} and have suppressed from the notation all canonical evaluation isomorphisms. The extra factor of $-1$ in the second equality is from the grading isomorphism $J$, which appears in $\tilde{\theta}$, and hence also $\tilde{\theta}_-$. On the other hand, since $\sigma_2 \sigma_1 \in G$, the right hand side of the desired equality is also the identity. Calculations for the remaining cases of $\pi(\sigma_1)$ and $\pi(\sigma_2)$ are similar.
\end{proof}
	
\begin{corollary} \label{cor:RealKnorrerPeriodicitydgCat}
The quasi-equivalence
\[
\mathcal{K} \circ \mathcal{K}:\MF(R,w) \rightarrow \MF(R\pser{u_1,v_1,u_2,v_2},w + u_1v_1 + u_2 v_2)
\]
admits a Real $G$-equivariant structure, where $\hat{G}$ is extended to act on $R\pser{u_1,v_1,u_2,v_2}$ by $\sigma(u_i)=\pi(\sigma) u_i$ and $\sigma(v_i)=v_i$, $\sigma \in \hat{G}$, and both categories are viewed as Real $2$-representations of $G$ by $(\rho, \theta)$.
\end{corollary}

\begin{proof}
This follows by applying \cref{thm:oriKnorrer} twice and noting that the universal twisting $\pi^*\hat{c} \in Z^2(\hat{G}; k^{\times}_{\pi})$ has order two.
\end{proof}

\begin{corollary}\label{cor:extendedKnorrerPeriodicityCat}
There is dg form quasi-equivalence
\[
\Perf(\MF_G(R,w))
\rightarrow
\Perf(\MF_G(R\pser{u_1,v_1,u_2,v_2},w + u_1v_1 + u_2 v_2),
\]
where both dg categories are given the dg duality structure of \cref{cor:dgDualityPerf}.
\end{corollary}

\begin{proof}
Applying \cref{thm:equivFunctorContra} to the quasi-equivalence $\mathcal{K} \circ \mathcal{K}$ of \cref{cor:RealKnorrerPeriodicitydgCat} gives a dg form functor
\[
((\mathcal{K} \circ \mathcal{K})^G,\varphi): \MF_G(R,w) \rightarrow \MF_G(R\pser{u_1,v_1,u_2,v_2},w + u_1v_1 + u_2 v_2))
\]
whose second component is a natural isomorphism. In view of \cref{ex:Real2RepVsDuality}, \cref{prop:Real2RepOnPerf} applied to the terminal $\Ctwo$-graded group gives a dg form functor
\[
(\Ind (\mathcal{K} \circ \mathcal{K})^G,\Ind \varphi): \Perf(\MF_G(R,w))
\rightarrow
\Perf(\MF_G(R\pser{u_1,v_1,u_2,v_2},w + u_1v_1 + u_2 v_2))
\]
whose first component is a quasi-equivalence by \cref{prop:fixedPointFunctorQusiEquiv} and \cref{rem:goodChar}.
\end{proof}

If $w$ is an isolated hypersurface singularity at the origin, \cref{cor:extendedKnorrerPeriodicityCat} applied to the terminal $\Ctwo$-graded group implies an equivalence
\[
\HMF(R,w) \simeq \HMF(R\pser{u_1,v_1,u_2,v_2},w + u_1v_1 + u_2 v_2)
\]
of triangulated categories with duality. This is precisely the extended Kn\"{o}rrer periodicity of Hori--Walcher \cite[\S 4.5]{hori2008}. \cref{cor:extendedKnorrerPeriodicityCat} therefore provides a precise mathematical formulation and generalization of Hori--Walcher's extended Kn\"{o}rrer periodicity.

\cref{thm:oriKnorrer} and its corollaries generalize to the twisted setting of \cref{ex:LGModelDiscTor}. The key observation is that the coherence condition \eqref{eq:equivariantFunctorContra} in the twisted case follows from the untwisted case, since the twisted coherence isomorphisms $\theta_{\sigma_2,\sigma_1}$ and $\tilde{\theta}_{\sigma_2,\sigma_1}$ are the same scalar multiple of their untwisted counterparts, namely $\hat{\mu}([\sigma_2 \vert \sigma_1])$. For brevity, we state only the twisted analogue of \cref{cor:extendedKnorrerPeriodicityCat}.

\begin{corollary}\label{cor:extendedKnorrerPeriodicityTwisted}
Let $\hat{\mu} \in Z^2(\hat{G}; k^{\times}_{\pi})$ with restriction $\mu \in Z^2(G; k^{\times})$. There is a dg form quasi-equivalence
\[
\Perf(\MF_{G,\mu}(R,w))
\rightarrow
\Perf(\MF_{G,\mu}(R\pser{u_1,v_1,u_2,v_2},w + u_1v_1 + u_2 v_2)),
\]
where both dg categories are given the dg duality structure of \cref{thm:dualityEquivMF}, but with $\theta_{\sigma,\sigma}$ replaced with $\hat{\mu}([\sigma \vert \sigma]) \theta_{\sigma,\sigma}$.
\end{corollary}

In the context of \cref{ex:LGModelDiscTor}, \cref{cor:extendedKnorrerPeriodicityTwisted} is an equivalence of categories of $D$-branes, with orientifold data determined by $(\hat{G}, \hat{\mu})$, for the Landau--Ginzburg $G$-orbifolds of $(R,w)$ and $(R\pser{u_1,v_1,u_2,v_2},w + u_1v_1 + u_2 v_2)$, both with discrete torsion $\mu$.

\appendix
\section{Calculations for the proof of \cref{thm:dualityEquivCat}} \label{app:A}

This appendix contains verifications of claims made in the proof of \cref{thm:dualityEquivCat}, whose notation we keep. For simplicity, we sometimes denote $\rho(\sigma)$ by $\sigma$.

\begin{lemma}\label{lem:functorOnFixedPoints}
The pair $(\rho(\sigma)(C),\{v_g\}_{g \in G})$ defines an object of $\catC^G$.
\end{lemma}

\begin{proof}
We need to verify the equality \eqref{eq:equivariantObject}, which states
$v_{g_2g_1} = \theta_{g_2,g_1} \circ g_2(v_{g_1}) \circ v_{g_2}$ for each $g_1,g_2 \in G$. The right hand side is
\[
\theta_{g_2,g_1} \circ g_2(v_{g_1}) \circ v_{g_2}
=
\theta_{g_2,g_1} \circ g_2\left(\theta_{g_1,\sigma}^{-1} \circ \theta_{\sigma, \sigma^{-1} g_1 \sigma} \circ \sigma(u^{-1}_{\sigma^{-1} g_1 \sigma}) \right) \circ \theta_{g_2,\sigma}^{-1} \circ \theta_{\sigma, \sigma^{-1} g_2 \sigma} \circ \sigma(u^{-1}_{\sigma^{-1} g_2 \sigma}),
\]
which by equation \eqref{eq:theta2CocycleContra} and naturality of $\theta_{g_2,\sigma}$ is equal to
\[
\theta^{-1}_{g_2 g_1,\sigma} \circ \theta_{g_2,g_1\sigma} \circ g_2(\theta_{\sigma, \sigma^{-1} g_1 \sigma}) \circ \theta_{g_2,\sigma}^{-1} \circ g_2\sigma(u^{-1}_{\sigma^{-1} g_1 \sigma}) \circ \theta_{\sigma, \sigma^{-1} g_2 \sigma} \circ \sigma(u^{-1}_{\sigma^{-1} g_2 \sigma}).
\]
Equation \eqref{eq:equivariantObject} with $g_i \mapsto \sigma^{-1} g_i \sigma$, $i=1,2$, and naturality of $\theta_{\sigma, \sigma^{-1} g_2 \sigma}$ then gives
\[
\theta_{g_2 g_1,\sigma}^{-1} \circ \theta_{\sigma, \sigma^{-1} g_2 g_1 \sigma} \circ \sigma(\theta_{\sigma^{-1} g_2 \sigma,\sigma^{-1} g_1 \sigma}^{-1}) \circ \theta_{\sigma, \sigma^{-1} g_2 \sigma}^{-1} \circ g_2 \sigma(u^{-1}_{\sigma^{-1} g_1 \sigma}) \circ \theta_{\sigma, \sigma^{-1} g_2 \sigma} \circ \sigma(u^{-1}_{\sigma^{-1} g_2 \sigma}).
\]
By repeated application of equation \eqref{eq:theta2CocycleContra}, this is equal to $v_{g_2g_1}$.
\end{proof}

\begin{lemma}\label{lem:ThetaMorphInCG}
The map $\Theta(\sigma)$ is a morphism in $\catC^G$.
\end{lemma}

\begin{proof}
The statement amounts to the equality \eqref{eq:morphismofhomotopyfixedpoints}, which reads $w_g \circ \Theta(\sigma) = g(\Theta(\sigma)) \circ u_g$ for each $g \in G$. We compute
\begin{eqnarray*}
w_g \circ \Theta(\sigma)
&=&
\theta_{g,\sigma}^{-1} \circ \theta_{\sigma, \sigma^{-1} g \sigma} \circ \sigma(\theta_{g,\sigma}) \circ \sigma(\theta_{\sigma, \sigma^{-1} g \sigma}^{-1}) \circ \sigma (\sigma(u_{\sigma^{-2} g \sigma^2})) \circ \theta_{\sigma,\sigma}^{-1} \circ u_{\sigma^2} \\
&=&
\theta_{g,\sigma}^{-1} \circ \theta_{\sigma, \sigma^{-1} g \sigma} \circ \sigma(\theta_{g,\sigma}) \circ \sigma(\theta_{\sigma, \sigma^{-1} g \sigma}^{-1}) \circ \theta_{\sigma,\sigma}^{-1} \circ \theta_{\sigma^2, \sigma^{-2} g \sigma^2} \circ u_{g \sigma^2} \\
&=&
g(\Theta(\sigma)) \circ u_g.
\end{eqnarray*}
The first equality follows from the definitions, the second from naturality of $\theta_{\sigma,\sigma}$ and euqation \eqref{eq:equivariantObject} with $g_1 \mapsto \sigma^{-2} g \sigma^2$ and $g_2 \mapsto \sigma^2$ and the third from equation \eqref{eq:theta2CocycleContra}.
\end{proof}

\begin{lemma}\label{lem:dgFormNatural}
The map $\varphi^{\sigma_1, \sigma_2}_{(C,u)} : \rho(\sigma_1)(C,u) \rightarrow \rho(\sigma_2)(C,u)$ is a morphism in $\catC^G$.
\end{lemma}

\begin{proof}
Write $\{v^{\sigma_1}_g\}_{g \in G}$ and $\{v^{\sigma_2}_g\}_{g \in G}$ for the homotopy fixed point data of $\sigma_1 C$ and $\sigma_2 C$, respectively. The present lemma amounts to the equality \eqref{eq:morphismofhomotopyfixedpoints}, which reads
$v_g^{\sigma_2} \circ \varphi^{\sigma_1, \sigma_2}_{(C,u)} = g(\varphi^{\sigma_1, \sigma_2}_{(C,u)}) \circ v_g^{\sigma_1}$ for each $g \in G$. The definitions give
\[
v_g^{\sigma_2} \circ \varphi^{\sigma_1, \sigma_2}_{(C,u)}
=
g(\sigma_2 (u_{\sigma_2^{-1} \sigma_1}) \circ \theta^{-1}_{\sigma_2, \sigma_2^{-1} \sigma_1} ) \circ \theta_{g,\sigma_1}^{-1} \circ \theta_{\sigma_1, \sigma_1^{-1} g \sigma_1} \circ \sigma_1 (u_{\sigma_1^{-1} g \sigma_1}^{-1}).
\]
Pre-applying $\sigma_2 \sigma_2^{-1}$ to $\sigma_1 (u_{\sigma_1^{-1} g \sigma_1}^{-1})$ gives
\begin{eqnarray*}
\sigma_1 (u_{\sigma_1^{-1} g \sigma_1}^{-1})
&=&
\theta_{\sigma_2,\sigma_2^{-1}} \circ \sigma_2 (\theta_{\sigma_2^{-1},\sigma_1} \circ \sigma_2^{-1} \sigma_1 (u_{\sigma_1^{-1} g \sigma_1}^{-1}) \circ \theta_{\sigma_2^{-1},\sigma_1}^{-1}) \circ \theta_{\sigma_2,\sigma_2^{-1}}^{-1}.
\end{eqnarray*}
Using equation \eqref{eq:equivariantObject} to rewrite $\sigma_2^{-1} \sigma_1 (u_{\sigma_1^{-1} g \sigma_1}^{-1})$ gives
\[
\sigma_2^{-1}\sigma_1 (u_{\sigma_1^{-1} g \sigma_1}^{-1}) = \theta_{\sigma_2,\sigma_2^{-1}} \circ \sigma_2 (\theta_{\sigma_2^{-1},\sigma_1} \circ u_{\sigma_2^{-1} \sigma_1} \circ u^{-1}_{\sigma_2^{-1} g \sigma_1} \circ \theta_{\sigma_2^{-1} \sigma_1, \sigma_1^{-1} g \sigma_1} \circ \theta_{\sigma_2^{-1},\sigma_1}^{-1}) \circ \theta_{\sigma_2,\sigma_2^{-1}}^{-1}.
\]
Similarly, we may pre-apply $\sigma_2 \sigma_2^{-1}$ to $g(\sigma_2 (u_{\sigma_2^{-1} \sigma_1}))$ to get a new expression for the latter.
Plugging these expressions into $v_g^{\sigma_2} \circ \varphi^{\sigma_1, \sigma_2}_{(C,u)}$ and repeatedly applying equation \eqref{eq:theta2CocycleContra} gives
\[
g(\varphi^{\sigma_1, \sigma_2}_{(C,u)}) \circ v_g^{\sigma_1}
=
\theta_{g,\sigma_2}^{-1} \circ \theta_{\sigma_2, \sigma_2^{-1} g \sigma_2} \circ \sigma_2 (u_{\sigma_2^{-1} g \sigma_2}^{-1}) \circ \sigma_2 (u_{\sigma_2^{-1} \sigma_1}) \circ \theta^{-1}_{\sigma_2, \sigma_2^{-1} \sigma_1},
\]
as required.
\end{proof}

\begin{lemma}\label{lem:dgFormCoher}
The pair $(\id_{\catC^G},\varphi^{\sigma_1,\sigma_2})$ is a dg form functor.
\end{lemma}

\begin{proof}
It remains to verify the the coherence condition \eqref{eq:formFuncCoher}, which states
\[
\sigma_2(\varphi^{\sigma_1, \sigma_2}_{(C,u)}) \circ \Theta(\sigma_2)_{(C,u)}
=
\varphi^{\sigma_1, \sigma_2}_{\sigma_1(C,u)} \circ \Theta(\sigma_1)_{(C,u)}
\]
for each $(C,u) \in \catC^G$. We have
\begin{eqnarray*}
\sigma_2(\varphi^{\sigma_1, \sigma_2}_{(C,u)}) \circ \Theta(\sigma_2)_{(C,u)}
&=&
\sigma_2(\sigma_2 (u_{\sigma_2^{-1} \sigma_1}) \circ \theta^{-1}_{\sigma_2, \sigma_2^{-1} \sigma_1}) \circ \theta_{\sigma_2,\sigma_2}^{-1} \circ u_{\sigma_2^2}\\
&=&
\sigma_2(\theta_{\sigma_2,\sigma_2^{-1} \sigma_1}^{-1}) \circ \theta_{\sigma_2,\sigma_2}^{-1} \circ \sigma_2^2(u_{\sigma_2^{-1} \sigma_1}) \circ u_{\sigma_2^2} \\
&=&
\sigma_2(\theta_{\sigma_2,\sigma_2^{-1} \sigma_1}^{-1}) \circ \theta_{\sigma_2,\sigma_2}^{-1} \circ \theta^{-1}_{\sigma_2^2, \sigma_2^{-1} \sigma_1} \circ u_{\sigma_2 \sigma_1}\\
&=&
\theta^{-1}_{\sigma_2,\sigma_1} \circ u_{\sigma_2 \sigma_1}.
\end{eqnarray*}
The first through fourth equalities follow by definition, naturality of $\theta_{\sigma_2,\sigma_2}$, equation \eqref{eq:equivariantObject} and equation \eqref{eq:theta2CocycleContra}, respectively. Similarly,
\begin{align*}
\varphi^{\sigma_1, \sigma_2}_{\sigma_1(C,u)} \circ \Theta(\sigma_1)_{(C,u)}
&=
\sigma_2 (v^{\sigma_1}_{\sigma_2^{-1} \sigma_1}) \circ \theta^{-1}_{\sigma_2, \sigma_2^{-1} \sigma_1} \circ \theta_{\sigma_1,\sigma_1}^{-1} \circ u_{\sigma_1^2} \\
&=
\begin{multlined}[t]
\theta^{-1}_{\sigma_2,\sigma_1} \circ \sigma_2 \sigma_1(u^{-1}_{\sigma_1^{-1} \sigma_2^{-1} \sigma_1^2}) \circ \theta_{\sigma_2,\sigma_1} \circ \sigma_2(\theta_{\sigma_1, \sigma_1^{-1} \sigma_2^{-1} \sigma_1^2}) \circ \sigma_2(\theta_{\sigma_2^{-1} \sigma_1,\sigma_1}^{-1}) \\
 \circ \theta^{-1}_{\sigma_2, \sigma_2^{-1} \sigma_1} \circ \theta_{\sigma_1,\sigma_1}^{-1} \circ u_{\sigma_1^2}
\end{multlined}\\
&=
\theta^{-1}_{\sigma_2,\sigma_1} \circ \sigma_2 \sigma_1(u^{-1}_{\sigma_1^{-1} \sigma_2^{-1} \sigma_1^2}) \circ \theta^{-1}_{\sigma_2 \sigma_1, \sigma_1^{-1} \sigma_2^{-1} \sigma_1^2}  \circ u_{\sigma_1^2}.
\end{align*}
The first, second and third equalities follow by definition, the explicit expression for $v^{\sigma_1}_{\sigma_2^{-1} \sigma_1}$ and naturality of $\theta_{\sigma_2,\sigma_1}$, and equation \eqref{eq:theta2CocycleContra}, respectively. Finally, equation \eqref{eq:equivariantObject} with $g_1 \mapsto \sigma_1^{-1} \sigma_2^{-1} \sigma_1^2$ and $g_2 \mapsto \sigma_2 \sigma_1$ gives $\varphi^{\sigma_1, \sigma_2}_{\sigma_1(C,u)} \circ \Theta(\sigma_1)_{(C,u)} = \theta^{-1}_{\sigma_2,\sigma_1} \circ u_{\sigma_2 \sigma_1}$.
\end{proof}

\begin{lemma}\label{lem:dgFormFunctTorsor}
For each $\sigma_1, \sigma_2, \sigma_3 \in \hat{G} \setminus G$, there are equalities of dg form functors
\[
(\id_{\catC^G}, \varphi^{\sigma_1, \sigma_3})
=
(\id_{\catC^G}, \varphi^{\sigma_2, \sigma_3}) \circ (\id_{\catC^G}, \varphi^{\sigma_1, \sigma_2}).
\]
\end{lemma}

\begin{proof}
We need to verify that $\varphi^{\sigma_1, \sigma_3}_{(C,u)}=\varphi^{\sigma_2,\sigma_3}_{(C,u)} \circ \varphi^{\sigma_1,\sigma_2}_{(C,u)}$ for each $(C,u) \in \catC^G$. We have
\begin{eqnarray*}
\varphi^{\sigma_1, \sigma_3}_{(C,u)}
&=&
\sigma_3(u_{\sigma_3^{-1} \sigma_1}) \circ \theta^{-1}_{\sigma_3, \sigma_3^{-1} \sigma_1} \\
&=&
\sigma_3(u_{\sigma_3^{-1}\sigma_2}) \circ \theta_{\sigma_3,\sigma_3^{-1} \sigma_2}^{-1} \circ \sigma_2 (u_{\sigma_2^{-1} \sigma_1}) \circ \theta_{\sigma_3,\sigma_3^{-1} \sigma_2} \circ \sigma_3 (\theta_{\sigma_3^{-1} \sigma_2, \sigma_2^{-1} \sigma_1}) \circ \theta^{-1}_{\sigma_3, \sigma_3^{-1} \sigma_1}.
\end{eqnarray*}
The first equality follows by definition and the second by equation \eqref{eq:equivariantObject} with $g_1 \mapsto \sigma_2^{-1} \sigma_1$ and $g_2 \mapsto \sigma_3^{-1} \sigma_2$. Applying equation \eqref{eq:theta2CocycleContra} to the fourth and fifth terms on the right hand side then gives $\varphi^{\sigma_2,\sigma_3}_{(C,u)} \circ \varphi^{\sigma_1,\sigma_2}_{(C,u)}$.
\end{proof}

\emergencystretch=1em
\bibliographystyle{amsalpha}

\newcommand{\etalchar}[1]{$^{#1}$}
\providecommand{\bysame}{\leavevmode\hbox to3em{\hrulefill}\thinspace}
\providecommand{\MR}{\relax\ifhmode\unskip\space\fi MR }
\providecommand{\MRhref}[2]{%
  \href{http://www.ams.org/mathscinet-getitem?mr=#1}{#2}
}
\providecommand{\href}[2]{#2}


\end{document}